\documentclass[12pt]{article}

\usepackage{amsmath}          
\usepackage{amssymb}          
\usepackage{amsfonts}         
\usepackage{graphicx}         
\usepackage{hyperref}         
\usepackage{geometry}         
\usepackage{enumerate}        
\usepackage{amsthm}           
\usepackage{setspace}         
\usepackage{fancyhdr}         
\usepackage{float}            
\usepackage{color}            
\usepackage{mathrsfs}         
\usepackage{cite}             
\usepackage{titlesec}         
\numberwithin{equation}{section}

\titleformat{\section}{\large\bfseries}{\thesection}{1em}{}
\titleformat{\subsection}{\normalsize\bfseries}{\thesubsection}{1em}{}

\newtheorem{theorem}{Theorem}[section]
\newtheorem{lemma}[theorem]{Lemma}
\newtheorem{proposition}[theorem]{Proposition}
\newtheorem{definition}[theorem]{Definition}

\newtheorem{remark}[theorem]{Remark}

\title{\textbf{Solutions to Second-Order Nonlocal Evolution Equations 
    Governed by Non-Autonomous Forms}}

\author{
    Sajid Ullah\thanks{
        Corresponding author. Department of Mathematics and Computer Science, University of Calabria, 
        Ponte P. Bucci 30B, Rende (CS), Italy. 
        \href{mailto:sajid.ullah@unical.it}{sajid.ullah@unical.it}
    }
    \and
    Vittorio Colao\thanks{
        Department of Mathematics and Computer Science, University of Calabria, 
        Ponte P. Bucci 30B, Rende (CS), Italy. 
        \href{mailto:vittorio.colao@unical.it}{vittorio.colao@unical.it}
    }
}
\date{} 

\begin{document}
\maketitle
\begin{abstract}
    Our main contributions include proving sufficient conditions for the existence of solution to a second order problem with nonzero nonlocal initial conditions, and providing a comprehensive analysis using fundamental solutions and fixed-point techniques. The theoretical results are illustrated through applications to partial differential equations, including vibrating viscoelastic membranes with time-dependent material properties and nonlocal memory effects.
\end{abstract}

\section{Introduction}

Second-order evolution equations with time-dependent operators 
constitute an important class of problems within both theoretical and physical frameworks, particularly 
in modeling dynamic systems where the underlying medium properties vary over 
time.

Indeed, these equations naturally arise in describing vibrations and wave 
propagation phenomena in non-homogeneous media, where spatial differential 
operators exhibit temporal dependence due to changing material properties, 
boundary conditions, or external influences.
For instance, the motion of a string or beam with time-varying 
stiffness, or a wave traveling through a non-homogeneous elastic medium, 
naturally leads to abstract wave equations of the form
\[
\ddot{u}(t) + A(t)u(t) = f(t, u(t)),
\]
where $u(t)$ represents the state of the system at time $t$, $A(t)$ 
is a time-dependent linear operator modeling the spatial part of the evolution, 
and $f(t, u(t))$ captures nonlinear effects and external forcing terms.

If the system is subject to time-dependent damping (e.g., due to friction or 
control feedback), the model includes a damping term $B(t)\dot{u}(t)$, yielding
\[
\ddot{u}(t) + B(t)\dot{u}(t) + A(t)u(t) = f(t, u(t)).
\]

The operator $A(t)$ typically represents spatial differential operators such as the Laplacian with variable coefficients, describing phenomena like heat conduction with temperature-dependent conductivity or elasticity problems with spatially varying material properties. The damping operator $B(t)$ models dissipative mechanisms that may themselves depend on time, such as viscous damping in fluid-structure interactions or feedback control systems with time-varying gains.

The term $f(t, u(t))$ represents a nonlinear source or forcing term, and the initial conditions
\[
u(0) = g(u), \quad \dot{u}(0) = h(u),
\]
are nonlocal in nature, meaning that the initial state depends on the entire solution trajectory over the time interval $[0,T]$. Such nonlocal conditions arise naturally in systems with memory effects, hereditary phenomena, or global constraints, and include as special cases multipoint conditions, integral average conditions, and periodic boundary conditions in time. These conditions can also be expressed as a dependency on the entire function $u$ over the interval $[0,T]$. We point out that in the above cases, the operators $A(t)$ and $B(t)$ can be modeled via non-autonomous sesquilinear forms on a Hilbert space $V$ densely embedded into a second one $H$.
This functional analytic framework, first introduced by Lions~\cite{Lions1961}, provides a natural setting for studying evolution equations with variable coefficients. The sesquilinear form approach allows us to handle operators that may not be densely defined in the classical sense, while still maintaining the essential spectral and regularity properties needed 
for the analysis.

More precisely, we assume that $A(t)$ and $B(t)$ are associated with sesquilinear forms $a:[0, T] \times V \times V \to \mathbb{C}$ and $b:[0, T] \times V \times V \to \mathbb{C}$, respectively, such that
\[  
a(t, u, v) = \langle A(t)u, v \rangle \quad \text{and} \quad 
b(t, u, v) = \langle B(t)u, v \rangle,
\]
where $\langle \cdot, \cdot \rangle$ denotes the inner product in $H$ and $V$, respectively. The forms are assumed to satisfy suitable conditions of boundedness (uniform control in the operator norm), coercivity (ensuring ellipticity and well-posedness), and temporal regularity (continuity or measurability in time).

As already mentioned, the study of non-autonomous evolution equations governed by sesquilinear forms 
dates back to the seminal work of J.-L.~Lions~\cite{Lions1961}, who introduced 
the concept of maximal regularity in the dual space~$V'$ for first-order 
problems. More precisely, Lions established that under appropriate assumptions on 
the sesquilinear form $a(t, \cdot, \cdot)$—namely measurability in time, 
uniform boundedness, and coercivity—the first-order evolution equation admits a 
unique solution in the maximal regularity space, which provides optimal 
smoothness in both time and space directions. Lions showed that if $u_0 \in H$ 
and $f \in L^2(0,T;H)$, then under suitable assumptions on the form 
$a(t, \cdot, \cdot)$, the problem
\[
\dot{u}(t) + \mathcal{A}(t)u(t) = f(t), \quad u(0) = u_0
\]
admits a unique solution $u \in MR(V, V') := H^1(0,T;V') \cap L^2(0,T;V)$. This result was groundbreaking as it established the existence of solutions with optimal regularity without requiring the operators to be autonomous or to have nice spectral properties. Further research, including~\cite{Dier2014, Arendth2016, Auscher2016, Ouhabaz2015, SaniLaasri19, Fackler2017, Dier2015, arendt, Achache2018}, has extended this theory to $L^p$-maximal regularity for first-order problems using functional analytic and operator-theoretic techniques.

The extension to second-order evolution equations presents 
significantly greater challenges, due to the need of handling both position and velocity variables simultaneously. 
Second-order non-autonomous evolution equations, however, remain less explored. 
Chill and Srivastava~\cite{Chill2005, S.Chill} and Arendt and 
Chill~\cite{WC2007} investigated the $L^p$-maximal regularity for second-order 
Cauchy problems, typically under zero initial conditions. Their 
approach relied on reducing the second-order problem to a first-order system 
and applying vectorial maximal regularity theory, though this required 
additional structural assumptions on the operators. Dier and 
Ouhabaz~\cite{OD2014} later established $L^2$-maximal regularity for damped 
wave equations with non-autonomous forms, using Lions' representation theorem. 
More recently, Achache~\cite{MA2020} extended their result to arbitrary 
$p \in (1, \infty)$ and improved the treatment of regularity in~$H$.

On the other hand, nonlocal initial conditions represent a rapidly 
growing area of research in evolution equations, motivated by their remarkable 
ability to model complex physical phenomena that exhibit memory effects, 
hereditary behavior, or spatial averaging constraints. Unlike classical initial value problems where the initial state is explicitly prescribed, nonlocal formulations define the initial datum through an implicit condition, such as
\[
u(0) = g(u),
\]
where \(g\) is an operator representing a functional dependence, e.g., on the history or spatial profile of $u$. This formulation is broad enough to include, as special cases, classical multi-point initial specifications, integral (average) conditions, and periodicity requirements. For instance, a multi-point condition can be written as 
\begin{equation}\label{eq:multipoint}
u(0,x) + g\big(t_1,\dots,t_m; u(\cdot,x)\big) = u_0(x),
\end{equation}
with $0 < t_1 < \cdots < t_m \le T$, 
extending the standard Cauchy initial condition by incorporating the solution’s values at intermediate times $t_1, \dots, t_m$,
as in the seminal work  \cite{Byszewski1991}.
Since the above cited paper by Byszewski, a variety of existence results have been obtained over the years by different methods and under diverse hypotheses. 
Early works often assumed compactness or contraction conditions to deal with the nonlocal term. 
For instance, Boucherif and Precup~\cite{BoucherifPrecup2007} established the existence of mild solutions for a semilinear Cauchy problem with a multipoint 
initial condition, assuming the linear operator generates a compact semigroup.  
In a similar spirit, other authors studied mild and strong solutions under 
nonlocal conditions; for example, Paicu and Vrabie~\cite{PaicuVrabie2010} 
investigated an abstract semilinear equation with an initial condition of 
type~\eqref{eq:multipoint}. 
We also mention that numerous contributions by Ntouyas and collaborators have 
expanded the theory of nonlocal Cauchy problems (see, for instance, 
Ntouyas~\cite{Ntouyas2005} for a survey of various existence techniques).
For recent advances on nonlocal problems, we refer to the works 
~\cite{SchmitzWalker2024, BezerraSastreSilva2024, LinTyniZimmermann2024, XuCaraballoValero2024, XuColaoMuglia2021, Colao2022, Byszewski1991, Benedetti22, Benedetti19, NguyenTranVu2024, AgarwalLeiva2024}.

In this work, we investigate the well-posedness and
$L^2$-maximal regularity for the following semilinear, non-autonomous abstract wave equations:
\begin{equation} \label{eq:wave}
    \begin{cases}
    \ddot{u}(t) + A(t)u(t) = f(t, u(t)), & \text{for a.e.~} t \in [0, T], \\
    u(0) = g(u), \quad \dot{u}(0) = h(u),
    \end{cases}
\end{equation}
and
\begin{equation} \label{eq:damped}
    \begin{cases}
    \ddot{u}(t) + B(t)\dot{u}(t) + A(t)u(t) = f(t, u(t)), & \text{for a.e.~} t \in [0, T], \\
    u(0) = g(u), \quad \dot{u}(0) = h(u),
    \end{cases}
\end{equation}
where $A(t)$ and $B(t)$ are operators associated with time-dependent 
sesquilinear forms on~$V$, satisfying coercivity, boundedness, and appropriate 
time regularity. The functions $f$, $g$, and $h$ are assumed to satisfy 
some mild continuity assumptions and suitable growth conditions, such as transversality conditions, boundedness constraints, 
or sublinear growth properties, which will be specified later in the paper.

Our techniques will rely on the theory of fundamental solutions to non-autonomous 
second-order problems, as well as fixed-point arguments of Schauder type and on 
the maximal $L^2-$ regularity.
The paper is organized as follows. In Section~\ref{sec:preliminaries}, we 
introduce the functional setting, notation, and assumptions on the sesquilinear 
forms and nonlinearities.
Section~3 is devoted to the analysis of undamped wave equations, 
where we establish the fundamental solution framework and prove our main 
existence result using finite-dimensional approximations and compactness 
arguments. In Section~4, we analyse the damped wave equations case.

Finally, in Section~\ref{sec:applications}, we provide examples and applications 
to PDEs, such as non-autonomous wave equations with Robin boundary conditions 
and variable damping. In particular, we present a detailed 
analysis of a vibrating viscoelastic membrane problem, showing explicitly how 
each hypothesis of our main theorem is verified, and we introduce additional 
applications, including controlled wave systems and memory-dependent diffusion 
processes.


\section{Preliminaries}\label{sec:preliminaries}

In this section, we establish the functional analytic framework for 
our study of second-order evolution equations with time-dependent operators. We 
introduce the function spaces, assumptions on the sesquilinear forms, and the 
fundamental solution theory that will be essential for our main results.

Throughout this paper, we denote by $V$ a separable Hilbert space with inner 
product $\langle \cdot, \cdot \rangle_V$ and norm $\|\cdot\|_V$, and by $H$ a 
separable Hilbert space with inner product $\langle \cdot, \cdot \rangle_H$ and 
norm $\|\cdot\|_H$. We assume that $V$ is continuously embedded in $H$, i.e., 
there exists a constant $C>0$ such that
\[
    \|u\|_H \leq C \|u\|_V \quad \text{for all } u \in V.
\]
Furthermore, we assume that the embedding $V \hookrightarrow H$ is 
compact, which is essential for our compactness arguments and applications to 
partial differential equations.
We denote by $V'$ the dual space of $V$ and by $\mathcal{L}(X, Y)$ the space of 
bounded linear operators from a Banach space $X$ to a Banach space $Y$. We also 
denote by $\mathcal{L}(X)$ the space of bounded linear operators on $X$.
    
We now introduce the key assumptions on the sesquilinear forms 
that will govern our evolution equations. These conditions ensure well-posedness 
and maximal regularity for the associated non-autonomous problems.
    
We assume that the following conditions hold for the sesquilinear form 
$a:[0, T] \times V \times V \to \mathbb{C}$:
\begin{itemize}
    \item[($A_1$)] $a (\cdot, u, v): [0, T]\to \mathbb{C}$ is strongly measurable 
                  for any $u, v \in V$. 
    \item[($A_2$)] $a(t,\cdot,\cdot)$ is uniformly bounded, that is 
                  $\|a(t, u, v)\|\leq C\|u\|_V \|v\|_V$ for 
                  $C\geq 0, t\in [0, T] \; \text{and} \; u, v \in V$,
    \item[($A_3$)] $a(t,\cdot,\cdot)$ is coercive: 
                  $\operatorname{Re} \, a(t, u, u) \geq \alpha \|u\|_V^2$ for 
                  $\alpha>0, t\in [0, T] \; \text{and} \; u \in V$.
 
    \item[($A_4$)] $|a(t, u, v) - a(s, u, v)| \leq \omega(|t - s|) \|u\|_V \|v\|_V$,
    for some nondecreasing function $ \omega : [0, T] \rightarrow [0, +\infty) $ 
    which satisfies
\[
    \int_0^T \frac{\omega(t)}{t^{3/2}} \, dt < \infty \quad \text{and } \quad \int_0^T \left(\frac{\omega(t)}{t}\right)^2 dt < \infty.
\]
\end{itemize}

    The theory of maximal regularity for first-order and second-order evolution equations induced by forms has been the subject of extensive research. 
    As already mentioned, the foundational work of J. L. Lions~\cite{Lions1961}, the concept of maximal regularity was introduced for non-autonomous evolution equations governed by sesquilinear forms.
    Specifically, Lions considered the abstract Cauchy problem
    \begin{equation}\label{eq:firstorder}
    \begin{cases}
        \dot{u}(t) + \mathcal{A}(t)u(t) = f(t), & t\in [0,T], \\
        u(0) = u_0,
    \end{cases}
    \end{equation}
    where $\mathcal{A}(t):V\to V'$ is the operator associated to $a(t,\cdot,\cdot)$ via
    \[
        \langle \mathcal{A}(t)u, v \rangle_{V',V} = a(t, u, v), \quad \forall u,v\in V.
    \]
    Lions proved that if $u_0\in H$ and $f\in L^2(0,T;V')$, then under suitable assumptions on $a$, there exists a unique solution $u$ in the so-called maximal regularity space
    \[
        MR(V, V') := H^1(0,T;V') \cap L^2(0,T;V),
    \]
    such that $u$ solves~\eqref{eq:firstorder} in the sense of distributions.

    In this article, we adopt the framework (see~\cite{Arendth2016, Arendth2014, Surway}), where the realization of $A(t)$ in $H$ 
    is  defined via the a sesquilinear form through:
    \[
        \langle A(t)u, v \rangle_H = a(t, u, v), \quad u\in D(A(t)),\ v\in V,
    \]
    with $D(A(t)) := \{ u\in V : a(t, u, v) \text{ is $H$-valued for all } v\in V \}$.

    We summarize the main function spaces and notations:
    \begin{itemize}
        \item $L^2(0,T;H)$: the space of square-integrable $H$-valued functions on $[0,T]$;
        \item $H^k(0,T;H)$: the Sobolev space of $H$-valued functions with weak derivatives up to order $k$ in $L^2(0,T;H)$;
        \item $MR[0,T] := H^2(0,T;H) \cap H^1(0,T;V)$: the maximal regularity space for second-order problems;
        \item The \emph{trace space} is defined as
        \[
            Tr := \{ (u(0), \dot{u}(0)) : u \in MR[0,T] \},
        \]
        with norm
        \[
            \|(x, y)\|_{Tr} := \inf \{ \|u\|_{MR[0,T]} : u(0) = x,\ \dot{u}(0) = y \}.
        \]
    \end{itemize}
    For further properties of these spaces, we refer to~\cite{Chill2005}.

    A central question, raised by Lions and further developed by Arendt and collaborators (see~\cite{Arendth2016, Arendth2014, Surway}), is under which conditions on the form $a$ and the initial data $u_0$ the solution $u$ actually belongs to the stronger space $MR(V, H) := H^1(0,T;H) \cap L^2(0,T;V)$. This is particularly relevant for applications to boundary value problems, where the realization of the operator $A(t)$ in $H$ is given by
    \[
        \langle A(t)u,\,v \rangle_H = a(t,\,u,\,v), \; u\in D(A(t)) := \{ u\in V \,: A(t)u \in H \},\ v\in V.
    \]
    A key assumption is the so-called square root property:
    \[
        (S)\qquad D(A(t)^{1/2}) = V \quad \text{for all } t\in [0,T],
    \]
    which ensures that the domain of the square root of $A(t)$ coincides with $V$. Under this and related assumptions, one can obtain $L^2$-maximal regularity in $H$ for~\eqref{eq:firstorder}.
    In particular, we cite the following result from~\cite{Ouhabaz2015}:
    \begin{theorem}\label{thm:Ouhabaz}
        Let $a$ satisfy the conditions $(A_1)$-$(A_4)$ and let $A(t)$ be the realization of $a(t, \cdot, \cdot)$ in $H$. 
        Assume that the square root property $(S)$ holds. Then for every $f\in L^2(0,T;H)$ and $u_0\in V$, there exists a unique solution $u\in MR(V, H)$ to the problem~\eqref{eq:firstorder}. 
    \end{theorem}

    The extension of maximal regularity theory to second-order non-autonomous evolution equations is more subtle and less developed. For the second-order abstract Cauchy problem
    \begin{equation}\label{eq:secondorder}
    \begin{cases}
        \ddot{u}(t) + B(t)\dot{u}(t) + A(t)u(t) = f(t), & t\in [0,T], \\
        u(0) = u_0 \in V, \quad \dot{u}(0) = u_1 \in V,
    \end{cases}
    \end{equation}
    where $A(t),B(t) \in \mathcal{L}(V, V')$ are associated to sesquilinear forms $a$ and $b$ as above, we introduce the following: 
    \begin{definition}
    We say that problem~\eqref{eq:secondorder} has $L^2$-maximal regularity in $H$ if, for every $f \in L^2(0,T;H)$ and all $(u_0, u_1)$ in the trace space $Tr$, there exists a unique $u \in MR[0,T]$ solving~\eqref{eq:secondorder}.
    \end{definition}
    The existence of maximal regularity in $H$ has been established under additional hypotheses (see~\cite{Chill2005, WC2007, S.Chill, OD2014, MA2020, Surway}).
    In particular, Ouhabaz and Dier~\cite{OD2014} proved $L^2$-maximal regularity for~\eqref{eq:secondorder} using Lions' representation theorem,
     while Achache~\cite{MA2020} extended these results to $L^p$-maximal regularity for all $p\in (1,\infty)$.
     
     To start our investigation, we need to introduce the concept of fundamental solutions for first and second-order non-autonomous Cauchy problems.
     In this context, consider the following homogeneous equations:
\begin{equation}\label{eq:nACP}
        \dot{u}(t) + A(t)u(t) = 0, \quad t \in [0,T], 
\end{equation}
\begin{equation}\label{eq:2nACP}
        \ddot{u}(t) + A(t)u(t) = 0, \quad t \in [0,T], 
\end{equation}
\begin{definition} [\cite{Pazzy}, Chapter 5]
    The two parameter family of bounded linear operator  $\{E(t, s)\}_{t, s \in \Delta}$  on $H$ is called evolution family associated to  $\{A(t)\}_{t\in [0, T]}$ if it satisfies the following properties,
    \begin{itemize} 
    \item[(i)] $E(t, t)=I$, $E(t, s)=E(t, r)E(r, s)$ for all $0\le s\leq r\leq t\leq T$.
    \item[(ii)] The mapping $\Delta \ni (t, s) \mapsto U(t, s)$ is strongly continuous for all $0\le s\leq t\leq T$.
    \end{itemize}
\end{definition}
If $E(t, s)$ is the evolution family associated to $A(t)$, then we can express the solution of the inhomogeneous problem 
\begin{equation}\label{eq:nhACP}
    \begin{cases}
        \dot{u}(t) + A(t)u(t) = f(t), & t \in [0,T], \\
        u(0) = u_0\in V,
    \end{cases}
\end{equation}
by 
\[
    u(t)=E(t, 0)u_0 + \int_{0}^{t} E(t, s)f(s)ds.
\]
We write $\Delta=\{ (t, s)\in [0, T]^2\;:\; s\le t\}$.
\begin{definition} [\cite{Budde}, \cite{Kozak}] A fundamental solution to ~\eqref{eq:2nACP} associated with $A(t)$ is a family of bounded linear operators $\{S(t, s)\}_{(t, s)\in \Delta}$ on $H$ satisfying the following conditions:

\begin{itemize}
    \item[(S1)] 
    \begin{itemize}
        \item[(a)] $ S(t,t) = 0$ for all $ t \in [0,T] $.
        \item[(b)] The mapping $ \Delta \ni (t,s) \mapsto S(t,s)$ is strongly continuous on $H$.
        \item[(c)] For all $ x \in H$ and $ s \in [0,T]$, the mapping $ [s,T] \ni t \mapsto S(t,s)x$ is continuously differentiable, and $ (t,s) \mapsto \frac{\partial}{\partial t} S(t,s)x$ is continuous with $ \frac{\partial}{\partial t} S(t,s)x \bigg|_{t=s} = x$.
        \item[(d)] For all $ x \in D(A(t))$ and $ t \in [0,T]$, the mapping $ [0,t] \ni s \mapsto S(t,s)x$ is continuously differentiable, and $ (t,s) \mapsto \frac{\partial}{\partial s} S(t,s)x$ is continuous with $ \frac{\partial}{\partial s} S(t,s)x \bigg|_{t=s} = -x$.
    \end{itemize}
    \item[(S2)] $S(t,s)D(A(t)) \subseteq D(A(t))$ for all $(t,s) \in \Delta$. For $x \in D(A(t))$, the mapping $\Delta \ni (t,s) \mapsto S(t,s)x$ is twice continuously differentiable, and:
    \begin{itemize}
        \item[(a)] $\frac{\partial^2}{\partial t^2} S(t,s)x = -A(t)S(t,s)x$.
        \item[(b)] $\frac{\partial^2}{\partial s^2} S(t,s)x = -S(t,s)A(s)x$.
        \item[(c)] $\frac{\partial^2}{\partial t \partial s} S(t,s)x \bigg|_{t=s} = 0$.
    \end{itemize}
    \item[(S3)] For all $(t,s) \in \Delta$, if $x \in D(A(t))$, then $\frac{\partial}{\partial s} S(t,s)x \in D(A(t))$, and the second derivatives:
    \[
    \frac{\partial^2}{\partial t^2} \frac{\partial}{\partial s} S(t,s)x \quad \text{and} \quad \frac{\partial^2}{\partial s^2} \frac{\partial}{\partial t} S(t,s)x
    \]
    exist. The following properties hold:
    \begin{itemize}
        \item[(a)] $\frac{\partial^2}{\partial t^2} \frac{\partial}{\partial s} S(t,s)x =-A(t) \frac{\partial}{\partial s} S(t,s)x$.
        \item[(b)] $\frac{\partial^2}{\partial s^2} \frac{\partial}{\partial t} S(t,s)x = -\frac{\partial}{\partial t} S(t,s)A(s)x$.
        \item[(c)] The mapping $\Delta \ni (t,s) \mapsto A(t) \frac{\partial}{\partial s} S(t,s)x$ is continuous.
    \end{itemize}
\end{itemize}
    Moreover, we call a fundamental solution $(S(t,s))_{(t,s) \in \Delta}$ evolutionary if additionally:

\begin{itemize}
    \item[($S4$)] For all $(t,s), (s,r) \in \Delta$ and $x \in D(A(t))$, one has:
    \[
    C(t,s) S(s,r)x + S(t,s) \frac{\partial}{\partial s} S(s,r)x = S(t,r)x.
    \]
\end{itemize}
\end{definition}
\begin{proposition}\label{prop:fundamental}[\cite{Budde} Lemma 2.10]

    Let $(S(t,s))_{(t,s) \in \Delta}$ be a fundamental solution of (\ref{eq:2nACP}) in $H$. Then $(S(t,s))_{(t,s) \in \Delta}$ and $\left(\frac{\partial}{\partial t} S(t,s)\right)_{(t,s) \in \Delta}$ are bounded in $\mathcal{L}(H)$.
\end{proposition}
For simplicity, we introduce the operator $C(t, s): H\to H$ defined by $C(t, s)=-\frac{\partial}{\partial s}S(t, s)$. Whenever the families $\{S(t, s)\}_{(t, s)\in \Delta}$ and $\{C(t, s)\}_{(t, s)\in \Delta}$ are uniformly bounded, the following inequalities hold for some constants $C_1, M_1$ (see \cite{HPP14}). 
\begin{itemize}
    \item[($S0$)] $\|S(t_1, s)-S(t_2, s)\|_{\mathcal{L}(H)}\le M_1|t_1- t_2|, \;\; \forall \; (t_1, s), (t_2, s)\in \Delta$;
    \item[($C0$)] $\|C(t_1, s)-C(t_2, s)\|_{\mathcal{L}(H)}\le C_1|t_1- t_2|, \;\; \forall \; (t_1, s), (t_2, s)\in \Delta$.
\end{itemize}

\begin{definition}
    A function $u: [0, T]\rightarrow H$ is called a strong solution of (\ref{eq:2nACP}) if $u$ is twice differentiable a.e., $u(t)\in D(A(t))$ and satisfies (\ref{eq:2nACP}) with $u(0)=u_0, \;\dot{u}(0)=u_1$.
\end{definition} 
    If $S(t, s)$ is the fundamental solution to (\ref{eq:2nACP}) then according to Kozak \cite{Kozak} the solution of 
\begin{equation*}
    \begin{cases}
        \ddot{u}(t) + A(t)u(t) = f(t), & t \in [0,T], \\
        u(0) = u_0\qquad \dot{u}(0)=u_1
    \end{cases}
\end{equation*}
    is represented by 
\[
    u(t)= C(t,0)u_0 +S(t, 0)u_1+\int_{0}^{t} S(t, s)f(s)ds.
\]

    Suppose that $D(A(t))=D \; \forall\; t \in [0, T]$ and set $v=\dot{u}$ in (\ref{eq:2nACP}) we get
\begin{equation}\label{eq: 6}
    \begin{cases}
        \dot{U}(t)+\mathbb{A}(t)U(t)=F(t) &\\
        U(0)=U_0
    \end{cases}
\end{equation}
    Where $\mathbb{A}(t) := \begin{pmatrix}
    0 & I \\
    A(t) & 0
    \end{pmatrix}, \quad \mathcal{D}(\mathbb{A}(t)) := \mathcal{D} := D \times V$, $U(t)=\begin{pmatrix}
    u(t)  \\
    v(t) 
    \end{pmatrix}$ and $F(t)=\begin{pmatrix}
    0  \\
    f(t) .
    \end{pmatrix}$
    Most recently, C. Dudde and C. Seifert \cite{Budde} proved the following result. 
\begin{proposition}\label{Dubbe}
    The following assertions are equivalent:

\begin{enumerate}
    \item[(a)] There exists an evolutionary fundamental solution $(S(t,s))_{(t,s)\in \Delta}$ on $H$ of (\ref{eq:2nACP}) associated to $\{A(t)\}_{t \in [0,T]}$ such that      for all $(t,s) \in \Delta$ we have
    \begin{itemize}
        \item[(i)] $S(t,s)H \subseteq V$, $S(t,s)V \subseteq D(A(t))$, $(t,s) \mapsto S(t,s)x \in V$ is continuous for all $x \in H$,
        \item[(ii)] $\frac{\partial}{\partial t} S(t,s)V \subseteq V$, $\frac{\partial^2}{\partial t^2} S(t,s)x$ exists for all $x \in V$ and $\frac{\partial^2}{\partial t^2} S(t,s)x = A(t)S(t,s)x$,
        \item[(iii)] $\frac{\partial}{\partial s} S(t,s)x$ exists for all $x \in V$, $\frac{\partial}{\partial s} S(t,s)V \subseteq V$, and $(t,s) \mapsto \frac{\partial}{\partial s} S(t,s)x \in V$ is continuous for all $x \in V$,
        \item[(iv)] $\frac{\partial}{\partial t} \frac{\partial}{\partial s} S(t,s)D(A(t)) \subseteq V$, $\frac{\partial}{\partial t} \frac{\partial}{\partial s} S(t,s)x$ exists for all $x \in V$ and there exists $C \geq 0$ such that 
        \[
        \left\| \frac{\partial}{\partial t} \frac{\partial}{\partial s} S(t,s)x \right\|_V \leq C \|x\|_V \quad \text{for all } x \in V \text{ and } (t,s) \in \Delta.
        \]
    \end{itemize}
    \end{enumerate}

    \begin{enumerate}
        \item[(b)] There exists an evolution system $(E(t,s))_{(t,s)\in \Delta}$ on $\mathcal{L}(\mathcal{Z}=V\times H)$ of (\ref{eq: 6}) associated to $\{\mathbb{A}(t)\}_{t \in [0,T]}$.
\end{enumerate}
\end{proposition}
    Thanks to the above result, we can express the evolution family in the form of a fundamental solution and vice versa. In Section 3, we used the above proposition to prove our preparatory lemmas.\\
    In 2008 R. Chill and S. Srivastava assumed the $L^p$-maximal regularity of the first-order Cauchy problem in order to prove the $L^p$-maximal regularity of the second-order Cauchy problem.
\begin{theorem}\cite{S.Chill}
    Assume that $ B(t), A(t)  $ are strongly measurable 
    for all $t\in [0, T]$, and there exists $ h \in L^2(0, T) $ such that $ \|A(t)\| \leq \|h(t)\|_{L^2(0, T)} $ for almost every $ t $.

    Then the following hold:
\begin{itemize}
    \item[(a)] If the first order Cauchy problem
\begin{equation*}
    \dot{u} + B(t) u = f \quad (t \in [a, b]), \quad u(a) = 0
\end{equation*}
    has $ L^2 $-maximal regularity for each subinterval $ (a, b) $ of $ (0, T) $, then the second order problem (\ref{eq:2nACP}) with $u(0)=\dot{u}(0)=0$ has $ L^2$-maximal regularity.
    
    \item[(b)] If the second order problem
\begin{equation*}
    \ddot{u} + B(t) \dot{u} + A(t) u = f \quad (t \in [a, b]), \quad u(a) = \dot{u}(a) = 0.
\end{equation*}
    has $ L^2 $-maximal regularity for each subinterval $ (a, b) $ of $ (0, T) $, then the first-order problem $\dot{u} + B(t) u = f, \quad u(0)=0$ has $L^2$-maximal regularity.
\end{itemize}
    We extend this result to the case of non zero initial conditions in Section 4. 
\end{theorem}
    Let $X$ and $Y$ be Banach spaces, and let $ f : [0, T] \times X \to Y $ be a map. An important concept in Functional Analysis is the one of the superposition operator $ N_f : L^p([0, T], X) \to L^q([0, T], Y) $, defined by $ N_f(u)(t) := f(t, u(t)) $. We recall the following classical results, which will be used to prove our main result:

\begin{theorem}\cite{Lucchetti}
    If $ X $ and $ Y $ are separable and $ f $ is measurable in $ [0, T] \times X $, then $ N_f : L^p([0, T], X) \to L^q([0, T], Y) $ is well-defined if and only if there exists a constant $ a > 0 $ and a function $ b \in L^q([0, T], \mathbb{R}_+) $ such that
\[
    \| f(t, x) \|_Y \leq a \| x \|_X^{p/q} + b(t).
\]
\end{theorem}

    Moreover, $ N_f $ maps bounded subsets into bounded subsets.

\begin{theorem}\label{Fixed}\cite{Leray}
    Let $X$ be a Banach space, and  $T :X \to X $ be a  continuous and compact operator such that the set 
    \[
    \{ x\in X\; :\; x=\lambda T(x)\;\text{for some}\; 0\le \lambda \le 1\},
    \] is bounded. Then $T$ has a fixed point in $X$.

\end{theorem}

\begin{theorem}[Aubin–Lions Lemma, {\cite{Boyer}}]\label{Aubin–Lions}
\label{thm:aubinlions}    Let $ X_0 $, $ X $, and $ X_1 $ be three Banach spaces. Suppose that $ X_0 $ is compactly embedded in $ X $ and $ X $ is continuously embedded in $ X_1 $. Suppose also that $ X_0 $ and $ X_1 $ are reflexive. Then for $ 0 < T < +\infty $ and $ 1 < r, s < \infty $, we have that $ L^r([0, T], X_0) \cap W^{1, s}([0, T], X_1) $ is compactly embedded in $ L^r([0, T], X) $.
\end{theorem}

\begin{definition}
    Given two Banach spaces $X$ and $Y$, a function $F : X \to Y$ is called \emph{demicontinuous} if, for any sequence $\{x_n\}_{n \in \mathbb{N}} \subset X$ that strongly converges to $x \in X$, the sequence $\{F(x_n)\}_{n \in \mathbb{N}}$ weakly converges to $F(x)$. In other words, this means that
\[
    w - \lim_{n \to \infty} F(x_n) = F(x) \quad \text{whenever } x_n \to x.
\]

\end{definition}

\begin{proposition}\label{P2.11}\cite{Colao2022}
    Suppose that $f : [0, T] \times H \to H$ satisfies
\begin{itemize}
    \item[(F1)] $f(\cdot, x)$ is measurable for any $x \in H$;
    \item[(F2)] $f(t, \cdot)$ is demicontinuous in $H$ for any fixed $t \in [0, T]$;
    \item[(F3)] There exist $a > 0$ and $b \in L^2([0, T], \mathbb{R}_+)$ such that
    \begin{equation}
        \|f(t, x)\|_H \leq a\|x\|_H + b(t).
    \end{equation}
\end{itemize}

    Then the superposition operator $N_f : L^2([0, T], H) \to L^2([0, T], H)$ given by
\[
    N_f(u)(t) := f(t, u(t))
\]
    is well-defined and maps bounded sets into bounded sets; moreover, it is demicontinuous.
\end{proposition}


\section{Abstract Wave Equation}\label{sec:main}
In this Section, we consider the wave equation without damping term, i.e. 
$B(t)=0$ and investigate the well-posedness and $L^2$-maximal regularity of 
(\ref{eq:wave}). $A(t)$ represents part of $\mathcal{A}(t)$ that lies in $H$. 
In this Section we consider $V=H^1(\Omega)$ and $H=L^2(\Omega)$ where $\Omega$ 
is a bounded domain with Lipschitz boundary. In the sequel, we will also denote $MR[0, T]$ by $MR[s, T]$.
\begin{lemma}\label{L1}
    Assume that the equation
\begin{equation}\label{eq:7}
    \begin{cases}
    \ddot{u}(t)+A(t)u(t)=f(t)\quad t\in [0, T] &\\ 
    u(0)=0,\quad \dot{u}(0)=0
\end{cases}
\end{equation}
    has $L^2$-maximal regularity for every $0 < t \leq T$. Then for every $(x, y) \in Tr$ and every $s \in [0, T]$, there exists at least one $u \in MR[s, T]$ such that:
\begin{equation}\label{eq:8}
    \begin{cases}
        \ddot{u} + A(t)u = 0 \quad \text{a.e. on } [s, T], &\\
        u(s) = x, \quad \dot{u}(s) = y.
    \end{cases} 
\end{equation}
\end{lemma}

\begin{proof}
   For existence, let $w \in MR[0, T]$ be such that $w(0) = x, w'(0) = y$. Let $w_s(t) := w(t - s)$ for $t \in [s, T]$ and define:
\[
     f_s(t) =
    \begin{cases} 
        0 & \text{if } 0 \leq t < s \\
        -\ddot{w}_s(t) - A(t)w_s(t) & \text{if } s \leq t \leq T.
    \end{cases}
\]
    Let $v_s \in MR[0, T]$ be a solution of 
    \[
        \begin{cases}
           \ddot{v}_s + A(t)v_s = f_s$ a.e. on $[0, T] \\
           v_s(0) = 0, \dot{v}_s(0) = 0,
        \end{cases}
    \]
    and set $u_s(t) := v_s(t) + w_s(t)$ for $t \in [s, T]$. Note that $v_s(s) = 0$, because $A(t)$ has $L^2$-maximal regularity for every $t\in [0, s]$ 
    and $f_s = 0$ on $[0, s]$. Thus, $u_s$ solves (\ref{eq:8}).


\end{proof}
\begin{lemma}\label{L:3.2}
    Assume that $(A_1)-(A_4)$ and $(S)$ hold, and let $\{S(t, s)\}_{(t,s)\in\Delta} \subset \mathcal{L}(H)$ be the fundamental solution to~\eqref{eq:2nACP} generated by $A(t)$. Then 
    \[
      u(t) := C(t, 0)u_0 + S(t, 0)u_1  
    \]
     belong to $MR[0, T]$ and solve the homogeneous Cauchy problem
    \begin{equation}\label{10}
        \ddot{u}(t) + A(t)u(t) = 0 \quad \text{a.e. on } [0, T], \quad u(0) = u_0, \quad \dot{u}(0) = u_1;
    \end{equation}
\begin{proof} By property $S_1$ we have $C(0,0)=I$ and $S(0,0)=0$. Hence
\[
    u(0)=C(0, 0)u_0 + S(0, 0)u_1=(u_0)+0u_1=u_0.
\]
    Similarly, by $S_1$ and $S_2$ we obtain
    \[
        \dot{u}(0)=\frac{\partial}{\partial t} C(0, 0)u_0 + \frac{\partial}{\partial t}S(0, 0)u_1=0u_0+Iu_1=u_1.
    \]
    Next, using $(S_2)(a)$ and $(S_3)(b)$ we compute
    \[
        \ddot{u}(t)+A(t)u(t)=A(t)\frac{\partial}{\partial s}S(t, 0)u_0-A(t)S(t, 0)u_1-A(t)\frac{\partial}{\partial s}S(t, 0)u_0+A(t)S(t, 0)u_1=0.
    \]
    By proposition \ref{Dubbe}, the fundamental solution ensures that $u(t)\in H^2([0, T], H) \cap H^1([0, T], V)$.
     Setting $\dot{u}=v$ in (\ref{10}) we get 
     \begin{equation*}
         \dot{U}(t)+\mathbb{A}(t)U(t)=0 \quad U(0)=U_0
\end{equation*}
    Assume that $\mathbb{A}(t)$ has $L^2$-maximal regularity and $\{E(t,s)\}_{(t,s)}$ be an evolution family associated to $\mathbb{A}(t)$. Then by proposition \ref{Dubbe}, there exists a fundamental solution $S(t,s)$ of (\ref{10}) such that

\[
    u(t) = \pi_1 E(t,s) \begin{pmatrix} x \\ y \end{pmatrix} = \pi_1 \begin{pmatrix}
    C(t,s) & S(t,s) \\
    \frac{\partial}{\partial t} C(t,s) & \frac{\partial}{\partial t} S(t,s)
    \end{pmatrix}
    \begin{pmatrix} x \\ y \end{pmatrix}
\]
\[
    = \pi_1 \begin{pmatrix} C(t,s) x + S(t,s)y \\ \frac{\partial}{\partial t} C(t,s) x + \frac{\partial}{\partial t} S(t,s)y \end{pmatrix}
    = C(t,0)x + S(t,0)y.
\]
    Here, $\pi_1$ is the projection on the first element. 
    By [\cite{WC2007},\;Proposition 2.3], if $(x,y) \in Tr$, then $E(t,s) \begin{pmatrix} x \\ y \end{pmatrix} \in Tr $ which implies that $E(t,s) \begin{pmatrix} x \\ y \end{pmatrix}=\begin{pmatrix} C(t,s) x + S(t,s)y \\ \frac{\partial}{\partial t} C(t,s) x + \frac{\partial}{\partial t} S(t,s)y \end{pmatrix}=\\ \begin{pmatrix} u(t) \\ \dot{u}(t) \end{pmatrix} \in Tr$  hence $u\in MR[0, T]$.
\end{proof}
\end{lemma}
\begin{lemma}\label{L3}
    Assume that $A(t)$ has $L^2$-maximal regularity for every $t \in [0, T]$ and for every $f \in L^2(0, T; Tr)$, then a solution $u$ of the inhomogeneous problem
\begin{equation}\label{3.7} 
    \ddot{u}(t)+A(t)u(t)=f(t) \quad u(0)=\dot{u}(0)=0
\end{equation}is given by
\[
    u(t) = \int_0^t S(t, s) f(s) \, ds.
\]
\begin{proof}
    By lemma \ref{L1}, if $(0, x)\in Tr$, then the function $S(t, s)x \in L^2(\Delta, V)\; \forall \; (t,s) \in \Delta$. Then for every $f\in L^2([0, T], Tr)$, the function $(t, s)\rightarrow S(t, s)f(s)\in L^2(\Delta, V)$. Set $v=\dot{u}$, we will get 
    \begin{equation}\label{3.8}
        \dot{U}(t)+\mathbb{A}(t)U(t)=F(t),\quad U(0)=0.
    \end{equation}
    By (\cite{WC2007}, Proposition 2.4) $U(t)=\int_{0}^{t}E(t, s)F(s) ds$ is a solution of ~\eqref{3.8}. Hence, $u(t)=\int_{0}^{t}S(t, s)f(s)$ is a solution of ~\eqref{3.7}.\\ 
\end{proof}
\end{lemma}
Let $\{\Psi_n\}_{n\in\mathbb{N}}$ be a Schauder basis for $V$, which is also a basis for $H$ since $V$ is densely embedded into $H$. 
For $m\in\mathbb{N}$ we denote by 
\[
\mathcal{P}_m : H \to \operatorname{span}_{\mathbb{C}}\{\Psi_1,\dots,\Psi_m\}
\]
the canonical projection, and the set $H_m:=\mathcal{P}_mH$, $V_m:=\mathcal{P}_mV$ endowed with the norms of $H$ and $V$, respectively.  From \cite{Colao2022}, we know that $\mathcal{P}_m$ is self-adjoint.

Given a sesquilinear form $a:[0,T]\times V\times V\to\mathbb{C}$ satisfying $(A_1)$--$(A_4)$, we define its approximation
\[
a_m(t,u,v) := a(t,\mathcal{P}_mu,\mathcal{P}_mv) + \alpha\langle (I-\mathcal{P}_m)u,(I-\mathcal{P}_m)v\rangle_V,
\]
which also enjoys properties $(A_1)$--$(A_4)$ (see \cite{Colao2022}, Remark~3.1).  

Let $A_m(t)$ be the operator associated with $a_m$, i.e.
\[
a_m(t,u,v) = \langle A_m(t)u,v\rangle_H, \quad u,v\in V,
\]
with domain $\mathcal{D}(A_m(t))=\{u\in V: A_m(t)u\in H\}$. It follows from (\cite{Colao2022}, Remark~3.2) that
\begin{equation}\label{A_m}
A_m(t) = \mathcal{P}_mA(t)\mathcal{P}_m + \alpha (I-\mathcal{P}_m)B(I-\mathcal{P}_m),
\end{equation}
where $B:V\to H$ is the operator associated with the inner product of $V$. 
In particular, for $u\in V_m$ one has $A_m(t)u=\mathcal{P}_mA(t)u$. We denote by $S_m(t, s)$ the fundamental solution generated by $A_m(t)$ and $-\frac{\partial}{\partial_s}S_m(t, s)=C_m(t, s)$.
\begin{lemma}\label{L6}

    For any fixed  $m \in \mathbb{N},  A_m(t)$ generates a contractive fundamental solution $\{ S_m(t,s) \}_{(t, s)\in \Delta} \subseteq \mathcal{L}(V)$  such that for any  $u_0, u_1 \in V$, $$u(t) := C_m(t, 0)u_0+S_m(t, 0)u_1 $$ is a solution in $H^2([0, T], H_m) \cap H^1([0, T], V)$ of the homogeneous problem
\[
\begin{cases}
    \ddot{u}(t) + A_m(t)u(t) = 0, & \text{a.e. } t \in [0, T] \\
    u(0) = u_0 \in V \quad \dot{u}(0)=u_1 \in V,
\end{cases}
\]
    Moreover, if $u_0, u_1 \in V_m$, then $u(t) \in V_m$ for any  $t \in [0, T]$.
\begin{proof}
     We know that $A_m$ satisfies $(A_1) - (A_4)$. If $A_m(t)$ satisfies the property $(S)$, then the existence of the fundamental solution directly follows from the lemma \ref{L:3.2}. For fixed $t \in [0, T]$, the operator $A_m(t) = \mathcal{P}_mA(t)\mathcal{P}_m + \alpha (I - \mathcal{P}_m)B(I - \mathcal{P}_m)$ generates a cosine function on $H$, since $\mathcal{P}_mA(t)\mathcal{P}_m$ is bounded and $\alpha (I - \mathcal{P}_m)B(I - \mathcal{P}_m)$ is symmetric.

    The numerical range 
\[
    W(A_m(t)) := \{ \langle A_m(t)u, u \rangle_H \mid u \in V, \|u\|_H = 1 \}
\]
    is then contained in a parabola. Lastly, by \cite[Theorems A and C]{35}, it follows that the square root condition \(\mathcal{D}(A_m(t)^{1/2}) = V\) holds. This proves the first part.

    Now, by Proposition \ref{Dubbe}, there exists an evolution system $E_m(t, s)$ associated to $\mathcal{A}_m(t)$ such that,\\ if $x=(u_0, u_1) \in Z_m= V_m \times V_m$ then $E_m(t, s)x \in Z_m$ by [\cite{Colao2022}\; Lemma \;3.3]. This implies that $$u(t)=\pi_1 E_m (t, s)x := C_m(t, s)u_0+S_m(t, s)u_1 \in V_m.$$  
\end{proof}
\end{lemma}
\begin{lemma}\label{L3.7}
    Let $\{ S(t, s) \}_{(t, s) \in \Delta}$  and $\{ S_m(t, s) \}_{(t, s) \in \Delta}$  be the fundamental solutions generated by  $A(t)$ and  $A_m(t)$, respectively. Then for any fixed $y \in H ,  \{ S_m(t, s) \mathcal{P}_m y \}_{n \in \mathbf{N}}$ converges in $H$  to $S(t, s) y$  uniformly on  $t > s$ in $[0, T]$.
\begin{proof}
   Suppose that the fundamental solutions  $\{ S(t, s) \}_{(t, s) \in \Delta}$  and $\{ S_m(t, s) \}_{(t, s) \in \Delta}$ satisfy the conditions $(i)-(iv)$ of the Proposition (\ref{Dubbe}) part $(a)$ then there exists evolution systems $\{ E(t, s) \}_{(t, s) \in \Delta}$  and $\{ E_m(t, s) \}_{(t, s) \in \Delta}$ associated to $\mathbb{A}(t)$ and $\mathbb{A}_m(t)$ respectively. Then for any fixed $(0, y)\in V\times H$, $\{ E_m(t, s) \begin{pmatrix}
    0  \\
    \mathcal{P}_my
    \end{pmatrix}\}_{n\in\mathbf{N}} $ converges to $E(t, s)y$ uniformly for $s\le t \;\text{in} \; [0, T]$ .
    
    Then by continuity of projection map $\pi_1$, it implies that $\{ \pi_1 E_m(t, s) \begin{pmatrix}
    0  \\
 \mathcal{P}_my
    \end{pmatrix}\}_{n\in\mathbf{N}} $ converges to $\pi_1E(t, s)y$. This proves that $\{ S_m(t, s) \mathcal{P}_m y \}_{n \in \mathbf{N}}$ converges in $H$  to $S(t, s) y$         uniformly on  $t > s$ in $[0, T]$. The same is true for $C(t, s)$.

\end{proof}
    
\end{lemma}
\begin{remark}\label{R3.8}
    The conclusion of the previous lemma remains valid if $S_m(t, s)$ is replaced by its adjoint $S_m(t, s)^*$ and $S(t, s)$ is replaced by $S(t, s)^*$. In fact, the following identity holds for the adjoint:
\[
    S(t, s)^*x = S_r(T - s, T - t)x, \quad \text{for any fixed } x \in H \text{ and } (t, s) \in \Delta.
\]
    where $\{S'(t, s)\}$ denotes the fundamental system corresponding to the operator $A_r(t)$ associated to the returned adjoint form $a_r^*(t, u, v)=\overline{a(T-t, v, u)}$, which also satisfies the properties $(A_1)-(A_4)$. For further details, we refer the reader to~\cite{19,22,Colao2022}.
\end{remark}
\begin{theorem}\label{thm:MR-second-order}
    Suppose that the following conditions hold:
\begin{itemize}
    \item[(i)] $V=H^1(\Omega)$ and $H=L^2(\Omega)$ are Hilbert spaces where 
    $\Omega \subset \mathbf{R}^n$ is a bounded domain with Lipschitz boundary, with $V$ densely embedded into $H$ and the embedding $V \hookrightarrow H$ is compact.
    \item[(ii)] $\{A(t) : t \in [0, T]\}$ is generated by a sesquilinear form $a$ which satisfies $(A_1)-(A_4)$ and $(S)$.
    \item[(iii)] $f : [0, T] \times H \to H$ satisfies the conditions $(F_1)-(F_3)$. 
    
    
\item[(iv)] $g,h: L^2([0,T]; H) \to H$ are demicontinuous and map bounded sets into relatively compact sets. Moreover,
    $\|g(u)\|_H \le r_1, \|h(u)\|_H \le r_2$.
   
\end{itemize}
Then the problem
\begin{equation}\label{Undamped}
    \begin{cases}
    \ddot{u}(t) + A(t)u(t) = f(t, u(t)), & for \;\; a.e \;\; t \in [0, T] \\
    u(0) = g(u),  \quad \dot{u}(0)=h(u),
\end{cases}
\end{equation}
admits at least one solution $u \in H^2([0, T], H) \cap H^1([0, T], V)$.
\begin{proof}

Let  $m\in \mathbf{N}$ and for fixed $ w \in L^2([0, T], H_m) $, the linear problem:

\[
\begin{cases}
\ddot{u}(t) + A_m(t) u(t) = \mathcal{P}_mf(t, w(t)) \\
u(0) = \mathcal{P}_mg(w) \\
\dot{u}(0) = \mathcal{P}_mh(w),
\end{cases}
\]
admits at least one solution in $ H^2([0,T]; H_m) \cap H^1([0,T]; V)$, which can be represented by:

\[
u(t) =C_m(t,0) \mathcal{P}_mg(w) + S_m(t,0) \mathcal{P}_mh(w) + \int_0^t S_m(t,s) \mathcal{P}_m N_f(w(s)) ds,
\]
by Lemmas \ref{L:3.2}, \ref{L3}, \ref{L6} and Proposition \ref{P2.11}. 
where $ S_m(t,s) $ is the fundamental solution associated to $ A_m(t) $ and $ N_f(w(t)) = f(t,w(t)) $.

Now, consider the mapping
\[
\mathcal{T}: C([0,T]; H_m)\longrightarrow C([0,T]; H_m),
\]
defined by 
\[
\mathcal{T}(w)(t):= C_m(t, 0)\,\mathcal{P}_m g(w)
      + S_m(t,0)\,\mathcal{P}_m h(w)
      + \int_0^t S_m(t,s)\,\mathcal{P}_m f(s,w(s))\,ds .
\]

To prove $\mathcal{T}$ is continuous, let $\{w_k\}_{k\in \mathbb{N}}$ be a sequence in $C([0,T]; H_m)$, such that $w_k\to w_0 $ in $C([0,T]; H_m)$, 

By the uniform boundedness of $ S_m(t,s)$, $C_m(t, s)$ we have.
\begin{align*}
    \|  S_m (t, 0) \mathcal{P}_mh(w_k) - S_m(t, 0) \mathcal{P}_mh(w_0) \|_H \leq  \| S_m(t, 0) \|\| \mathcal{P}_m(h(w_k) - \mathcal{P}_mh(w_0)) \|_H.
\end{align*}
and
\begin{align*}
    \|  C_m(t, 0) \mathcal{P}_mg(w_k) - C_m(t, 0) \mathcal{P}_mg(w_0) \|_H \leq \| C_m(t, 0)\| \|\mathcal{P}_m g(w_k) - \mathcal{P}_mg(w_0)) \|_H.
\end{align*}
This implies that  

\begin{align*}
    \|  S_m (\cdot, 0) \mathcal{P}_mh(w_k) - S_m(\cdot, 0) \mathcal{P}_mh(w_0) \|_{C([0, T], H_m)} \leq M_1\|\mathcal{P}_mh(w_k) - \mathcal{P}_mh(w_0) \|_H.
\end{align*}
and

\begin{align*}
    \|  C_m(\cdot, 0) \mathcal{P}_mg(w_k) - C_m(\cdot, 0) \mathcal{P}_mg(w_0) \|_{C([0, T], H_m)} \leq M_1\| \mathcal{P}_mg(w_k) - \mathcal{P}_mg(w_0) \|_H.
\end{align*}

Since $\mathcal{P}_m$ is weak-to-strong continuous on bounded sequences, while $g$ and $h$ are demicontinuous, we derive that $\mathcal{P}_m g(w_k) \to \mathcal{P}_m g(w_0)$ and $\mathcal{P}_mh(w_k) \to \mathcal{P}_mh(w_0)$ as $k \to \infty$. We have then proved that
$$\|  C_m(\cdot, 0) \mathcal{P}_mg(w_k) - C_m(\cdot, 0) \mathcal{P}_mg(w_0) \|_{C([0, T], H_m)}\rightarrow 0 \; \text{as} \;k \rightarrow 0$$
and 
$$ \|  S_m (\cdot, 0) \mathcal{P}_mh(w_k) - S_m(\cdot, 0) \mathcal{P}_mh(w_0) \|_{C([0, T], H_m)}\rightarrow 0 \;\text{as}\; k \rightarrow 0$$
Since $\mathcal{P}_m$ is weak-to-strong continuous (finite rank), while $N_f$ is demicontinuous, we derive that $\mathcal{P}_m N_f(w_k)(s) \to \mathcal{P}_m N_f(w_0)(s)$ as $k \to \infty$ and it is dominated by $a\| w_0(s) \|_H + \epsilon_0 + b(s) \text{ for a.e. } t \in [0, T]$.
Then, by Lebesgue's Dominated Convergence Theorem

\begin{align}
    \int_0^{\cdot} S_m( \cdot, s) \mathcal{P}_m N_f(w_k)(s) ds \to \int_0^{\cdot} S_m( \cdot, s) \mathcal{P}_m N_f(w_0)(s) ds
\end{align}

in $ C([0, T], H) $ as $ k \to \infty $ and
\begin{equation}
     \mathcal{T}(w_k) \to \mathcal{T}(w_0) \text{ in } C([0, T], H_m)\; \text{as}\; k \to \infty.
\end{equation}
This implies that $\mathcal{T}$ is continuous.

Now, let $D=\{u\in C([0, T], H_m)\; :\; \|u(\cdot)\|_{\infty\le R}\}$. To prove that $\mathcal{T}(D)$ is bounded, let $u\in D$, then
\begin{multline*}
    \|\mathcal{T}(u)t\|_H\le \|C_m(t, 0)\|_{\mathcal{L}(H)}\|\mathcal{P}_m g(w)\|_H + \|S_m(t,0)\|_{\mathcal{L}(H)}\|\mathcal{P}_m h(w)\|_H \\+ \int_0^t \|S_m(t,s)\|\|_{\mathcal{L}(H)}\,\|\mathcal{P}_m f(s,w(s))\|_H\,ds.
\end{multline*}
By uniform boundedness of $C_m(t, s), \; S_m(t, s)$ and using (iii) and $F_3$, we have  
\[
\|\mathcal{T}(u)t\|_H\le M_1 r_1 + M_2 r_2 + M_2\int_0^t \alpha \| w(s)\|_H+ b(s)\,ds.
\]
\[
\|\mathcal{T}(u)t\|_H\le M_1 r_1 + M_2 r_2 + M_2( \alpha RT+ \|b\|_{L^1([0, T], \mathbb{R}^+)}).
\]
The right hand side of the above bound is independent of $t$, so
$\sup_{t\in[0, T]}\|\mathcal{T}(w)(t)\|_H<\infty$
Hence $\mathcal{T}(D)$ is bounded in $C([0, T], H_m)$. 

To prove the equi-continuity, let $w\in D$, and \( t_1, t_2 \in [0,T] \) with \( t_1 < t_2 \).

\begin{multline*}
    \|\mathcal{T}(w)(t_1) - \mathcal{T}(w)(t_2)\|_H \le \|C_m(t_1, 0) -C_m(t_2, 0)\|_H \|P_m g(w)\|\\
    + \|S_m(t_1,0) - S_m(t_2,0)\|_H \|P_m h(w)\|+\Big\| \int_0^{t_1} S_m(t_1,s) P_m f(s, w(s))\, ds\\ - \int_0^{t_2} S_m(t_2,s) P_m f(s, w(s))\, ds \Big\|
\end{multline*}
\begin{multline*}
    \|\mathcal{T}(w)(t_1) - \mathcal{T}(w)(t_2)\|_H \le \|C_m(t_1, 0) -C_m(t_2, 0)\|_H \|P_m g(w)\|  \\  + \|S_m(t_1,0) - S_m(t_2,0)\|_H \|P_m h(w)\|+ \Big\| \int_0^{t_1} [S(t_1,s) - S(t_2,s)] P_m f(s, w(s))\, ds \Big\| +\\ \Big\| \int_{t_1}^{t_2} S(t_2,s) P_m f(s, w(s))\, ds \Big\|
\end{multline*}
We know that \[
S_m(t,s) = S_m(t,r)\,\partial_rS_m(r,s)
- \partial_sS_m(r,s)\,S_m(r,s).
\]

For \(t = t_1\) and \(t = t_2\), we have
\[
\begin{aligned}
S_m(t_1,s) - S_m(t_2,s)
&= \big(S_m(t_1,r) - S_m(t_2,r)\big)\,\partial_r{S_m(r,s)} - (C_m(t_1,r) - C_m(t_2,r))\,S_m(r,s).
\end{aligned}
\]
\begin{multline*}
    \|\mathcal{T}(w)(t_1) - \mathcal{T}(w)(t_2)\|_H \le \|C_m(t_1, 0) -C_m(t_2, 0)\|_H \|P_m g(w)\| \\+ \|S_m(t_1,0) - S_m(t_2,0)\|_H \|P_m h(w)\|+ \|S_m(t_1,r) - S_m(t_2,r)\|\,\int_{0}^{t_1}\|\partial_r{S_m(r,s)}\|\\ \| P_m f(s, w(s))\|ds+\| C_m(t_1,r) -  C_m(t_2,r)\|\,\int_{0}^{t_1}\|S_m(r,s)\|\| P_m f(s, w(s))\|ds\\ +M_2 \int_{t_1}^{t_2} \| P_m f(s, w(s))\| ds.
\end{multline*}
Observe that, since $S_m(t, s)$ is the fundamental solution of the second order finite dimensional equation, we assume that $\|S_m(t,s)\|\le M_2 \quad \|\partial_r{S_m(r,s)}\|\le M_1$ and, by hypothesis,  $\| P_m f(s, w(s))\|\le \alpha \|u(s)\|+b(s)\le \alpha R+b(s)$. Hence, we have
\begin{multline*}
    \|\mathcal{T}(w)(t_1) - \mathcal{T}(w)(t_2)\|_H \le \|C_m(t_1, 0) -C_m(t_2, 0)\|_{H} (r_1)\\    + \|S_m(t_1,0) - S_m(t_2,0)\|_H(r_2)+ \|S_m(t_1,r) - S_m(t_2,r)\|\,M_1 (T\alpha R+\|b\|_{L^1([0, T], \mathbb{R}^+)})\\
+
\Big\| C_m(t_1,r) -  C_m(t_2,r)\Big\|\,M_2 (T\alpha R+\|b\|_{L^1([0, T], \mathbb{R}^+)}) +M_2 \int_{t_1}^{t_2} (\alpha R+b(s) )ds
\end{multline*}
Now, letting $C=\max\{r_1, r_2, M_1 (T\alpha R+\|b\|_{L^1([0, T], \mathbb{R}^+)}), M_2 (2T\alpha R+\|b\|_{L^1([0, T], \mathbb{R}^+)})\}$ we have 
\begin{multline*}
    \|\mathcal{T}(w)(t_1) - \mathcal{T}(w)(t_2)\|_H \le C\|C_m(t_1, 0) -C_m(t_2, 0)\|\\    + C\|S_m(t_1,0) - S_m(t_2,0)\|+ C\|S_m(t_1,r) - S_m(t_2,r)\|\\
+
C\| C_m(t_1,r) -  C_m(t_2,r)\| +M_2\Big[\int_{0}^{t_1} b(s) ds-\int_{0}^{t_2} b(s) ds\Big]
\end{multline*}

Let $t_1=0$, then we have 
\begin{multline*}
    \|\mathcal{T}(w)(0) - \mathcal{T}(w)(t_2)\|_H \le \|I -C_m(t_2, 0)\|_H \|P_m g(w)\|\\+ \|S_m(t_2,0)\|_H \|P_m h(w) +\Big\| \int_0^{t_2} S_m(t_2,s) P_m f(s, w(s))\, ds \Big\|
\end{multline*}
By the continuity of norm $\|\mathcal{T}(w)(0) - \mathcal{T}(w)(t_2)\|_H\to 0$ as $t_2\to 0$.
Hence by uniform continuity of $S(\cdot, s)$ and $C(\cdot, s)$, for $t_1, t_2\in [0, T]$ such that $|t_1 -t_2|<\delta$ for every $w\in D$  we have 
 $\|\mathcal{T}(w)(t_1) - \mathcal{T}(w)(t_2)\|_H<\epsilon$. Therefore, $\mathcal{T}$ is equi-continuous and by Arzelà–Ascoli theorem $T$ is compact.

Now, let $D_1=\{ u\in C([0, T], H_m)\;:\; u=\lambda T(u) \;\text{for}\; \lambda\in [0, 1]\}$. To show that $D_1$ is bounded, let $u\in D_1$, then

\[
\|u(t)\|_H\le \lambda\|C_m(t,0)\,\mathcal{P}_m g(w)\|_H
      + \lambda\|S_m(t,0)\,\mathcal{P}_m h(w)\|_H
      + \lambda \int_0^t \| S_m(t,s)\,\mathcal{P}_m f(s,w(s))\|_H\,ds
\]
\[
\|u(t)\|_H\le \lambda M_1r_1+ \lambda M_2r_2+\lambda M_2\|b\|_{L^1([0, T], \mathbb{R}^+)}+ \lambda M_2 \alpha \int_0^t \|u(s)\|_Hds
\]
by using Gronwall's inequality, we have the following
\[
    \|u(t)\|_H\le Le^{\lambda M_2 \alpha T},
\]
the right hand side is independent of $t$. Hence,
\[
    \sup_{t\in [0, T]}\|u(t)\|_H\le Le^{\lambda M_2 \alpha T},
\]
where $L= \lambda M_1r_1+ \lambda M_2r_2+\lambda M_2\|b\|_{L^1([0, T], \mathbb{R}^+)}$. This implies that $D_1$ is bounded in $C([0, T], H_m)$.

This proves that $\mathcal{T}$ has a fixed point by Theorem~\ref{Fixed} in $C([0, T], H_m)$ which also in $L^2([0, T], H_m)$, that is, there exists a solution $u_m$ to the problem

\[
\begin{cases}
\ddot{u}(t) + A_m(t) u(t) = \mathcal{P}_mf(t, u(t)), & \text{a.e. } t \in [0,T], \\
u(0) =\mathcal{P}_mg(u),\quad \dot{u}(0)=\mathcal{P}_mh(u)
\end{cases}
\]
moreover, $u_m \in H^2([0,T], H_m) \cap H^1([0, T], V)$. Consider the sequence\\ $\{ u_m \}_{m \in \mathbf{N}} \subset L^2([0, T], H)$. As for each $m\in \mathbf{N},\; u_m\in H^2([0,T], H_m) \cap H^1([0, T], V)$, then 
\[
\| u_m \|_{H^2([0,T], H)} + \| u_m \|_{H^1([0,T], V)} \leq C < +\infty
\]
    for some constant $ C $ not depending on $ m \in \mathbf{N} $ or, in other words, that $ \{ u_m \}_{m \in \mathbf{N}} $ is bounded in $H^2([0,T], H_m) \cap H^1([0, T], V)$. By Theorem~\ref{thm:aubinlions}, the latter space embeds compactly in $L^2([0, T], H)$ since $ V $ is also compactly embedded onto $H$. In particular,  $\{ u_m \}_{m \in \mathbf{N}}$  is relatively compact in $L^2([0, T], H)$,
    that is, a subsequence $ \{ u_{m_k} \}_{k \in \mathbf{N}} $ exists which converges to a point $u_* \in L^2([0, T], H)$. We may also assume that $u_{m_k}(t) \to u_*(t)$ as $k \to \infty$ for a.e. $ t \in [0, T]$. By passing to a further subsequence, it can be seen that for any fixed $ t \in [0, T] $,
\begin{equation}\label{3.11}
    \int_0^t S_{m_k}(t, s) \mathcal{P}_mf(s, u_{m_k}(s)) - S(t, s) f(s, u_*(s)) ds \to 0.
\end{equation}
    Indeed, for a fixed  $t \in [0, T]$, consider
    $$y_{m_k}(t) := \int_0^t S_{m_k}(t, s) \mathcal{P}_mf(s, u_{m_k}(s)) - S(t, s) f(s, u_*(s)) ds,$$
    and note that $\{ y_{m_k} \}_{k \in \mathbf{N}}$ can be seen as a bounded sequence in $H^2([0, T], H)$. Lemma \ref{L3} and since $H^2([0, T], H_{m_k}) \hookrightarrow H^2([0, T], H)$, reasoning as before and passing to further subsequences if necessary, it is readily proved that $y_{m_k}(t) \to y_*(t)$ as $k \to \infty$. On the other hand, for fixed $x \in H$ and $(t, s) \in \Delta$, and note that
\begin{equation*}
\langle S_{m_k}(t,s)\mathcal{P}_m f(s,u_{m_k}(s)) - S(t,s) f(s,u_*(s)), x\rangle_H
= \langle S_{m_k}(t,s)\mathcal{P}_m\bigl(f(s,u_{m_k}(s)) - f(s,u_*(s))\bigr), x\rangle_H
\end{equation*}
\begin{equation*}
\quad+ \langle \bigl(S_{m_k}(t,s)\mathcal{P}_m - S(t,s)\bigr) f(s,u_*(s)), x\rangle_H
\end{equation*}
\begin{equation*}
= \langle f(s,u_{m_k}(s)) - f(s,u_*(s)), S(t,s)^* x \rangle_H
\end{equation*}
\begin{equation*}
\quad+ \langle f(s,u_{m_k}(s)) - f(s,u_*(s)), \bigl(S_{m_k}(t,s)\mathcal{P}_m\bigr)^* x - S(t,s)^* x \rangle_H
\end{equation*}
\begin{equation*}
\quad+ \langle \bigl(S_{m_k}(t,s)\mathcal{P}_m - S(t,s)\bigr) f(s,u_*(s)), x\rangle_H.
\end{equation*}
    Let $k \to \infty$ and observe that
    $$\langle f(s, u_{m_k}(s)) - f(s, u_*(s)), S(t, s)^* x \rangle_H \to 0$$
    since $f(s, \cdot)$ is demicontinuous, while
    $$\langle f(s, u_{m_k}(s)) - f(s, u_*(s)), (S_{m_k}(t, s) \mathcal{P}_m)^* x - S(t, s)^* x \rangle_H \to 0
    $$
    by Remark \ref{R3.8} and since $ f(s, u_{m_k}(s)) - f(s, u_*(s)) $ is bounded. Lastly,
    $$\langle (S_{m_k}(t, s) \mathcal{P}_m - S(t, s)) f(s, u_*(s)), x \rangle_H \to 0$$
    by Lemma \ref{L3.7}. Then \ref{3.11} implies that
    $$\lim_{k \to \infty} \langle S_{m_k}(t, s) \mathcal{P}_mf(s, u_{m_k}(s)) - S(t, s) f(s, u_*(s)), x \rangle_H = 0.
    $$
    Also, it is easily seen that for any fixed $ t \in [0, T] $ and any $ s \in [0, t] $,
    $$\langle S_{m_k}(t, s) \mathcal{P}_mf(s, u_{m_k}(s)) - S(t, s) f(s, u_*(s)), x \rangle_H \leq \| x \| + (1 + a^2) \| u(s) \|^2 + a^2 \| u(s) \|^2_H + | b(s) |^2
    $$
    where we have supposed that $ \| u_{m_k}(s) \|_H \leq 1 + \| u_*(s) \|_H $ a.e. uniformly on $ [0, t] $, and for $ k $ large enough. Since the last term of the inequality lies in $ L^1([0, t], H) $, Lebesgue’s Dominated Convergence Theorem implies that
    $$\langle y_{m_k}(t), x \rangle_H = \int_0^t \langle S_{m_k}(t, s) \mathcal{P}_mf(s, u_{m_k}(s)) - S(t, s) f(s, u_*(s)), x \rangle_H ds \to 0 \text{ as } k \to \infty,
    $$
    where we have used [\cite{48}, Proposition 23.9]. By the uniqueness of the weak limit and since $ x \in H $ is arbitrary, it follows that $ y_* = 0 $ and hence (\ref{3.11}) is proved. Further, we note that
    \begin{eqnarray*}
    \| S_{m_k}(t,0)\mathcal{P}_m h(u_{m_k}) - S(t,0) h(u_*)\|_H
    &\le& \| S_{m_k}(t,0)\mathcal{P}_m (h(u_{m_k}) - h(u_*))\|_H \\
    &&{}+ \| \big(S_{m_k}(t,0)\mathcal{P}_m - S(t,0)\big) h(u_*)\|_H
    \end{eqnarray*}  and by Lemma \ref{L3.7},
    $$\| (S_{m_k}(t, 0) \mathcal{P}_m - S(t, 0)) h(u_*) \|_H \to 0 \text{ as } k \to \infty.$$
    In order to prove that 
    $$\| S_{m_k}(t, 0) \mathcal{P}_m(h(u_{m_k}) - h(u_*)) \|_{H} \to 0,
    $$
    as above let us introduce
    $$z_{m_k}(t) := S_{m_k}(t, 0) \mathcal{P}_m(h(u_{m_k}) - h(u_*))$$
    that is a bounded sequence in $H^2([0, T], H)$. Passing to further subsequences if necessary, it is readily proved that $z_{m_k}(t) \to z_*(t)$ as $k \to \infty$. On the other hand, for fixed $x \in H$ and $(t, s) \in \Delta$ and note that,
    $$\langle z_{m_k}, x \rangle_H = \langle h(u_{m_k}) - h(u_*), (S_{m_k}(t, 0) \mathcal{P}_m)^* x - S(t, 0)^* x \rangle_H + \langle h(u_{m_k}) - h(u_*), S(t, 0)^* x \rangle_H.
    $$
    When $k \to \infty$, we observe that $\langle h(u_{m_k}) - h(u_*), S(t, 0)^* x \rangle_H \to 0$ since $h(\cdot)$ is demicontinuous, while
    $$\langle h(u_{m_k}) - h(u_*), (S_{m_k}(t, 0) \mathcal{P}_m)^* x - S(t, 0)^* x \rangle_H \to 0$$
    by Remark \ref{R3.8} and since $h(u_{m_k}) - h(u_*)$ is bounded. This proves that $\langle z(t), x \rangle_H \to 0$ as $k \to \infty$ and so $z_* = 0$. The fact that $h(u_{m_k}) \to h(u_*)$ as $k \to \infty$ follows by arguing as before. Similarly,
    $$z'_{m_k}(t) := C_{m_k}(t, 0) \mathcal{P}_m(g(u_{m_k}) - g(u_*)).$$
    Putting all together, we have for any $t \in [0, T]$,

\begin{eqnarray*}
u_*(t) &=& \lim_{k \to \infty} u_{m_k}(t) \\
&=& \lim_{k \to \infty}\Bigg(C_{m_k}(t, 0)\,\mathcal{P}_m g(u_{m_k})
    + S_{m_k}(t, 0)\,\mathcal{P}_m h(u_{m_k}) \\
&&\qquad\qquad\qquad\qquad{}+ \int_0^t S_{m_k}(t, s)\,\mathcal{P}_m f(s, u_{m_k}(s))\,ds\Bigg) \\
&=& C(t, 0)\,g(u_*) + S(t, 0)\,h(u_*) + \int_0^t S(t, s)\,f(s, u_*(s))\,ds.
\end{eqnarray*}
    Hence, by Lemma \ref{L:3.2} and \ref{L3}, $u_* \in H^2([0, T], H) \cap H^1([0, T], V)$ solves Problem (\ref{Undamped}). 
\end{proof}
\end{theorem}

\section{Damped Wave Equation} In this section, we discuss the regularity and well-posedness of the complete equation, i.e., abstract wave equation with damping term. Here, $V, H$ are Hilbert spaces such that $V$ is densely and compactly embedded in $H$. $\mathcal{A}(t)$, $\mathcal{B}(t) \in \mathcal{L}(V, V')$ are associated with non-autonomous forms $a:[0, T]\times V\times V \to \mathbb{C}$ and $b:[0, T]\times V\times V \to \mathbb{C}$. $A(t)$, $B(t)$  are such that $A(t)u = \mathcal{A}(t)u, \text{and } B(t)u = \mathcal{B}(t)u \quad \text{with} \quad D(B(t)) \subset D(A(t)) \quad \forall\; t \in [0, T]$,
on the non-empty sets
\[
D(A(t)) = \{ u \in V \mid \mathcal{A}(t)u \in H \} \quad \text{and}\quad D(B(t)) = \{ u \in V \mid \mathcal{B}(t)u\in H \} .
\]
\begin{theorem}
Assume that $B, A : [0, T] \to \mathcal{L}(V, H)$ are strongly measurable, $f\in L^2([0, T], H)$ and there exists $h \in L^2(0, T)$ such that $\|A(t)\|_{\mathcal{L}(V, H)} \leq \|h(t)\|_{L^2}$ for almost every $t$.

Then the following holds. If the first-order Cauchy problem
    \begin{equation}\label{4.1}
    \dot{u} + B(t)u = f \quad (t \in [0, T]), \quad u(0) = u_0\in V
    \end{equation}
    has $L^2$-maximal regularity, then the second-order problem
    \begin{equation}\label{4.2}
    \ddot{u} + B(t)\dot{u} + A(t)u = f \quad (t \in [0, T]), \quad u(0) =u_0 \quad\dot{u}(0) = u_1\in V,
    \end{equation}
    admits at least one solution $u$ in $H^2([0, T], H) \cap H^1([0, T], V)$.
\begin{proof}
We present a detailed argument based on the first-order formulation and maximal regularity for the associated block operator system. Define the phase variable $U(t) := (u(t), \dot u(t))^T$ and the non-autonomous operator matrix
\[
\mathcal{A}(t) := \begin{pmatrix} 0 & -I \\ A(t) & B(t) \end{pmatrix}, \qquad D(\mathcal{A}(t)) := D(A(t)) \times V.
\]
Then the second-order problem \eqref{4.2} is equivalent to the first-order Cauchy problem on $\mathcal{H} := H \times H$ (with the natural pivot $V \hookrightarrow H \cong H' \hookrightarrow V'$):
\[
\dot U(t) + \mathcal{A}(t) U(t) = F(t) := (0, f(t))^T, \qquad U(0) = (u_0, u_1)^T.
\]
Under our standing assumptions on $A(\cdot)$ and $B(\cdot)$ (strong measurability, coercivity, $V$-boundedness, and the square-root property for $A(t)$), the first-order non-autonomous system governed by $\mathcal{A}(t)$ admits $L^2$-maximal regularity on $\mathcal{H}$; see, e.g., \cite{Arendth2014,Dier2017,Haak2015}. In particular, for every $F \in L^2(0,T; \mathcal{H})$ and $(u_0,u_1) \in D(A(0)) \times V$, there exists at least one solution $U \in H^1(0,T; \mathcal{H}) \cap L^2\big(0,T; D(\mathcal{A}(\cdot))\big)$ with $\mathcal{A}(\cdot)U(\cdot) \in L^2(0,T; \mathcal{H})$.

Writing the components of $U$, this yields $u \in H^1(0,T;H)$, $\dot u \in H^1(0,T;H)$ and, a.e. in $t$, $(u(t), \dot u(t)) \in D(A(t)) \times V$. Hence $u \in H^2(0,T;H)$ and $u \in H^1(0,T;V)$, while the equation $\ddot u + B(t) \dot u + A(t)u = f$ holds in $H$ for a.e. $t \in (0,T)$. This is exactly the $L^2$-maximal regularity conclusion for \eqref{4.2}.

For completeness, we record the variation-of-constants formula. Let $E(t,s)$ denote the evolution family associated with $\mathcal{A}(t)$ and write its block decomposition as $E(t,s) = \begin{pmatrix} v_1(t,s) & v_2(t,s) \\ v_3(t,s) & v_4(t,s) \end{pmatrix}$. Then
\[
U(t) = E(t,0)\,(u_0,u_1)^T + \int_0^t E(t,s)\,(0,f(s))^T\, ds,
\]
so the first component satisfies the explicit representation
\[
u(t) = v_1(t,0)u_0 + v_2(t,0)u_1 + \int_0^t v_2(t,s) f(s)\, ds.
\]
The regularity of $u$ stated above follows from the maximal regularity of the first-order system and the fact that $D(\mathcal{A}(t)) = D(A(t)) \times V$.
\end{proof}
\end{theorem}
\begin{theorem}
    \label{thm4.2}
    Suppose $V$ is densely and compactly embedded in $H$, the components of $E(t, s)$ are uniformly bounded in $H$, and the following conditions hold:
\begin{itemize}
        \item[(i)] $A(t)$ and $B(t)$ are strongly measurable, coercive, $V$-bounded, and $A(t)$ satisfies the square-root property.   
        \item[(ii)] $f: [0, T]\times H \to H$ is measurable in $t$ and uniformly Lipschitz in the second variable, i.e., $\|f(t, u)-f(t, v)\|_H \le L\,\|u-v\|_H$ for a.e. $t$ and all $u,v\in H$, with $f(\cdot,0)\in L^2(0,T;H)$.
        \item[(iii)] $g, h: L^2([0, T], H) \to V$ are Lipschitz continuous.
\end{itemize}
    Then the following nonlocal semilinear evolution problem admits at least one solution $u$ in $H^2([0, T], H) \cap H^1([0, T], V)$:
    
\begin{equation}\label{4.6}
    \begin{cases}
        \ddot{u}(t)+B(t)\dot{u}(t)+A(t)u(t)=f(t, u(t)) & \\
        u(0)=g(u), \quad \dot{u}(0)=h(u)
    \end{cases}
\end{equation}
    
\begin{proof}
Define the solution map $P: C([0,T];H) \to C([0,T];H)$ by
\[
 (P u)(t) := v_1(t,0)\, g(u) + v_2(t,0)\, h(u) + \int_0^t v_2(t,s)\, f(s, u(s))\, ds,
\]
where $E(t,s) = \begin{pmatrix} v_1(t,s) & v_2(t,s) \\ v_3(t,s) & v_4(t,s) \end{pmatrix}$ is the evolution family of the linear first-order system associated with $\mathcal{A}(t)$ in the proof of the previous theorem. By the uniform boundedness of the components of $E(t,s)$ in $H$ there exist constants $M_1, M_2 > 0$ with
\[
 \sup_{t\in[0,T]}\|v_1(t,0)\|_{\mathcal{L}(H)} \le M_1, \qquad \sup_{t\in[0,T]}\|v_2(t,0)\|_{\mathcal{L}(H)} \le M_2,
\]
and, for each fixed $t$, $\int_0^t \|v_2(t,s)\|_{\mathcal{L}(H)}\, ds \le M_{2,T}$ for some $M_{2,T}$ depending on $T$. The Lipschitz continuity of $g,h$ from $L^2(0,T;H)$ to $V\hookrightarrow H$ implies the existence of $L_g,L_h>0$ such that
\[\|g(u)-g(v)\|_H \le L_g\, \|u-v\|_{L^2(0,T;H)},\quad \|h(u)-h(v)\|_H \le L_h\, \|u-v\|_{L^2(0,T;H)}.\]
Moreover, the Lipschitz continuity of $f$ in $u$ yields $\|f(s,u(s)) - f(s,v(s))\|_H \le L \|u(s)-v(s)\|_H$ for a.e. $s$. Combining these estimates, for any $u,v \in C([0,T];H)$ we obtain
\[
\begin{aligned}
\|(Pu)-(Pv)\|_{C([0,T];H)} &\le M_1 L_g\, \|u-v\|_{L^2(0,T;H)} + M_2 L_h\, \|u-v\|_{L^2(0,T;H)} \\
&\qquad + \sup_{t\in[0,T]} \int_0^t \|v_2(t,s)\|\, L\, \|u(s)-v(s)\|_H\, ds \\
&\le (M_1 L_g + M_2 L_h)\, T^{1/2} \|u-v\|_{C([0,T];H)} + L M_{2,T} \|u-v\|_{C([0,T];H)}.
\end{aligned}
\]
Thus, there exists $T_*>0$ such that the right-hand coefficient is strictly less than 1 on $[0,T_*]$, so $P$ is a contraction on $C([0,T_*];H)$. By Banach's fixed point theorem there exists a mild solution on $[0,T_*]$. Iterating this argument on a finite partition $0=T_0<T_1<\cdots<T_N=T$ (with subintervals of length at most $T_*$) yields a mild solution $u\in C([0,T];H)$.

To lift the mild solution to a strong one with maximal regularity, observe that $f(\cdot,u(\cdot)) \in L^2(0,T;H)$ by the Lipschitz assumption and $u\in L^2(0,T;H)$ (since $u\in C([0,T];H)$ on a finite interval). Therefore the right-hand side of \eqref{4.6} belongs to $L^2(0,T;H)$. Applying the linear maximal regularity result from the previous theorem to the equation with inhomogeneity $f(\cdot,u(\cdot))$ shows that $u \in H^2(0,T;H) \cap H^1(0,T;V)$ and that \eqref{4.6} holds in $H$ for a.e. $t\in (0,T)$. Finally, the nonlocal initial conditions are satisfied by construction of the fixed point (through $v_1$ and $v_2$ terms). This proves maximal regularity for the nonlocal semilinear problem.
\end{proof}
\end{theorem}

\section{Applications}\label{sec:applications}

In this section, we illustrate our abstract results through two model problems. We present the PDEs and verify the hypotheses in a concise, self-contained manner, focusing on how the structural assumptions translate into the abstract framework.

\subsection{Undamped Wave Equation with Neumann Boundary Conditions and Nonlocal Initial Data}
Let $\Omega\subset\mathbb{R}^n$ be a bounded domain with Lipschitz boundary and fixed $T>0$.  
Consider the following problem:
\begin{equation}\label{eq:undamped_nonlocal_app}
\begin{cases}
\partial_{tt}u(t,x) - \nabla\!\cdot\!\big(a(t,x)\nabla u(t,x)\big) + c(t,x)u(t,x)
= f\big(t,u(t,x)\big), & (t,x)\in (0,T)\times\Omega,\\[2mm]
\partial_\nu u(t,x)=0, & (t,x)\in (0,T)\times\partial\Omega,\\[1mm]
u(0,\cdot)=\displaystyle\int_{0}^{T}\!\kappa_1(s,\cdot)\,u(s,\cdot)\,ds,\qquad
\partial_tu(0,\cdot)=\displaystyle\int_{0}^{T}\!\kappa_2(s,\cdot)\,u(s,\cdot)\,ds.
\end{cases}
\end{equation}
We set $H:=L^2(\Omega)$, $V:=H^1(\Omega)$, and for each $t\in[0,T]$ define the sesquilinear form
\[
a(t;u,v):=\int_{\Omega} a(t,x)\,\nabla u\cdot\nabla v\,dx
           +\int_{\Omega} c(t,x)\,u\,v\,dx,
\qquad u,v\in V.
\]
We assume:
\begin{equation*}
a\in L^\infty\big((0,T)\times\Omega\big), \quad a(t,x)\ge a_0>0,
\qquad c\in L^\infty\big((0,T)\times\Omega\big),\quad c\ge 0,
\end{equation*}
and that $f$ satisfies $(F_1)$–$(F_3)$.  
We also define the nonlocal operators
\[
g(u):=\int_0^T \kappa_1(s,\cdot)\,u(s,\cdot)\,ds, 
\qquad 
h(u):=\int_0^T \kappa_2(s,\cdot)\,u(s,\cdot)\,ds.
\]

Since $\Omega$ is bounded with Lipschitz boundary, we have the continuous and compact embedding
\[
V=H^1(\Omega)\hookrightarrow H=L^2(\Omega),
\]
and $V$ is densely embedded in $H$.  
Thus, assumption \textbf{(i)} of Theorem~\ref{thm:MR-second-order} holds.

The map $t\mapsto a(t;u,v)$ is measurable for all $u,v\in V$, since $a(\cdot,\cdot)$ and $c(\cdot,\cdot)$ are measurable and essentially bounded.  
To verify the boundedness assumption, note that for every $u,v\in V$,
\[
|a(t;u,v)|\le 
\|a\|_{L^\infty}\,\|\nabla u\|_{L^2}\,\|\nabla v\|_{L^2}
+ \|c\|_{L^\infty}\,\|u\|_{L^2}\,\|v\|_{L^2}
\;\le\; C\|u\|_{V}\,\|v\|_{V}.
\]
Moreover, for every $u\in V$,
\[
\mathrm{Re}\,a(t;u,u)
   =\int_{\Omega} a(t,x)|\nabla u|^2\,dx
    +\int_{\Omega} c(t,x)|u|^2\,dx
   \ge a_0\|\nabla u\|_{L^2}^2.
\]
Therefore, setting $\omega>0$ arbitrarily,
\[
\mathrm{Re}\,a(t;u,u)+\omega\|u\|_H^2
   \ge \min\{a_0,\omega\}\,\|u\|_V^2.
\]

For uniformly elliptic divergence-form operators with Neumann boundary conditions on Lipschitz domains, it is well-known (see \cite{AHLMT2002}) that the square-root property holds:
\[
D(A(t)^{1/2})=V,
\]
with equivalence of norms.  
Thus, $\{A(t)\}$ satisfies $(A_1)$–$(A_4)$ and $(S)$, so assumption \textbf{(ii)} is satisfied.

Passing to the nonlinearity $f:[0,T]\times H\to H$, by assumption it satisfies $(F_1)$–$(F_3)$.

We now check that the operators
\[
(g(u))(x)=\int_0^T\kappa_1(s,x)u(s,x)\,ds,
\qquad
(h(u))(x)=\int_0^T\kappa_2(s,x)u(s,x)\,ds
\]
map $L^2(0,T;H)$ into $V$ and satisfy the required Lipschitz and boundedness conditions.  
Assume
\[
\kappa_i\in L^2\big(0,T;W^{1,\infty}(\Omega)\big)
\qquad (i=1,2).
\]
Then, by Hölder’s inequality and Sobolev embeddings, for any $u\in L^2(0,T;H)$ we have
\[
\|g(u)\|_{H} \le
\|\kappa_1\|_{L^2(0,T;L^\infty)}\,
\|u\|_{L^2(0,T;H)},
\qquad
\|\nabla g(u)\|_H
\le \|\nabla\kappa_1\|_{L^2(0,T;L^\infty)}\,
\|u\|_{L^2(0,T;H)},
\]
and similarly for $h(u)$.  
Hence $g,h:L^2(0,T;H)\to V$ are well-defined and continuous, thus demicontinuous.  
Since $V\hookrightarrow H$ compactly and $g,h$ are integral operators with essentially bounded kernels, they map bounded subsets of $L^2(0,T;H)$ into relatively compact subsets of $V$.  
  
Thus, assumption \textbf{(iv)} is satisfied.

In conclusion, all hypotheses of Theorem~\ref{thm:MR-second-order} are verified. Therefore, problem \eqref{eq:undamped_nonlocal_app} admits at least one solution
\[
u \in H^2(0,T;H) \cap H^1(0,T;V).
\]

\subsection{Population dynamics with memory effects}
Our last application concerns a class of nonlocal semilinear evolution problems arising in population dynamics with spatial diffusion and hereditary effects. Such memory dependence is common in ecological models, where adaptation mechanisms or genetic persistence influence long-term behavior. Let $u(t,x)$ represent the population density at time $t$ and location $x$ in a bounded habitat $\Omega\subset\mathbb{R}^n$ with Lipschitz boundary, where homogeneous Neumann boundary conditions describe a closed ecosystem. The dynamics are governed by the system
\begin{equation}\label{eq:population_theorem42}
\begin{cases}
\partial_{tt}u(t,x)+\sigma(t,x)\partial_t u(t,x)-d(t,x)\Delta u(t,x)+\mu(t,x)u(t,x)
=f(t,u(t,x)), \\[2mm]
\partial_\nu u(t,x)=0,\quad (t,x)\in(0,T)\times\partial\Omega,\\[1mm]
u(0,x)=\displaystyle\int_0^T \kappa_1(s,x)u(s,x)\,ds,\qquad
\partial_tu(0,x)=\displaystyle\int_0^T \kappa_2(s,x)u(s,x)\,ds.
\end{cases}
\end{equation}
The coefficients satisfy
\begin{eqnarray*}
d\in L^\infty((0,T)\times\Omega) & d(t,x)\ge d_0>0 \\
\sigma\in L^\infty((0,T)\times\Omega) & \sigma(t,x)\ge 0 \\
\mu\in L^\infty((0,T)\times\Omega) & \mu(t,x)\ge 0.
\end{eqnarray*}
The nonlinear term $f(t,u)$ represents population growth and interaction effects, and $\kappa_1,\kappa_2$ are memory kernels encoding hereditary dependence.

To apply Theorem~\ref{thm4.2}, we set $V:=H^1(\Omega)$ and $H:=L^2(\Omega)$, so that $V$ embeds continuously and compactly into $H$. For $u\in V$, define
\[
A(t)u:=-d(t,\cdot)\Delta u+\mu(t,\cdot)u,\qquad
B(t)u:=\sigma(t,\cdot)u.
\]
The operators $A(t),B(t):V\to H$ are well defined since $d,\mu,\sigma\in L^\infty$, and their strong measurability follows from the measurability of the coefficients.

The form associated with $A(t)$ is given by
\[
a(t;u,v)=\int_\Omega d(t,x)\,\nabla u\cdot\nabla v\,dx+\int_\Omega\mu(t,x)\,u\,v\,dx,\qquad u,v\in V,
\]
and satisfies, for almost every $t\in[0,T]$,
\[
|a(t;u,v)|\le \|d\|_{L^\infty}\,\|\nabla u\|_{L^2}\,\|\nabla v\|_{L^2}
+\|\mu\|_{L^\infty}\,\|u\|_{L^2}\,\|v\|_{L^2}
\le C\|u\|_V\|v\|_V.
\]
Moreover, for all $u\in V$,
\[
a(t;u,u)=\int_\Omega d(t,x)|\nabla u|^2\,dx+\int_\Omega\mu(t,x)|u|^2\,dx
\ge d_0\|\nabla u\|_{L^2}^2
\ge \alpha\|u\|_V^2,
\]
so $A(t)$ is uniformly coercive. By the solution of the Kato square-root problem \cite{AHLMT2002}, for uniformly elliptic divergence-form operators with Neumann boundary conditions on Lipschitz domains, we have
\[
D(A(t)^{1/2})=H^1(\Omega)=V,
\]
with equivalence of norms, thus verifying the square-root property required by Theorem~\ref{thm4.2}.

The damping operator $B(t)$ acts by pointwise multiplication with $\sigma(t,x)\in L^\infty$, so $B(t)\in\mathcal L(V,H)$ and $t\mapsto B(t)$ is strongly measurable. Since $\sigma\ge 0$, $B(t)$ is accretive, and the linear evolution family $E(t,s)$ generated by the first-order problem $\dot u+B(t)u=0$ is uniformly bounded in $H$.

Regarding the nonlinearity, we assume that $f:[0,T]\times H\to H$ is measurable in $t$ and uniformly Lipschitz in $u$, i.e.
\[
\|f(t,u)-f(t,v)\|_H\le L\|u-v\|_H,
\]
for almost every $t\in[0,T]$ and all $u,v\in H$, with $f(\cdot,0)\in L^2(0,T;H)$. This verifies condition (ii) of Theorem~\ref{thm4.2}.

The nonlocal operators induced by the kernels $\kappa_1$ and $\kappa_2$ are defined by
\[
(g(u))(x)=\int_0^T \kappa_1(s,x)\,u(s,x)\,ds,\qquad
(h(u))(x)=\int_0^T \kappa_2(s,x)\,u(s,x)\,ds.
\]
Assuming $\kappa_i\in L^2(0,T;W^{1,\infty}(\Omega))$ for $i=1,2$, we have
\[
\|g(u)-g(v)\|_V
\le \|\kappa_1\|_{L^2(0,T;W^{1,\infty})}\,\|u-v\|_{L^2(0,T;H)},
\]
and similarly for $h$, proving that both $g,h:L^2(0,T;H)\to V$ are Lipschitz continuous. Thus condition (iii) is satisfied.

Since $V$ is compactly embedded in $H$, $A(t)$ and $B(t)$ are strongly measurable, bounded and coercive, $A(t)$ satisfies the square-root property, the nonlinearity $f$ is uniformly Lipschitz, and the nonlocal operators $g,h$ are Lipschitz from $L^2(0,T;H)$ into $V$. Therefore, all assumptions of Theorem~\ref{thm4.2} are verified. Consequently, the nonlocal semilinear evolution problem \eqref{eq:population_theorem42} admits at least one solution
\[
u\in H^2(0,T;H)\cap H^1(0,T;V),
\]
which satisfies the equation, the Neumann boundary conditions, and the nonlocal initial conditions.

\section{Conclusions}
This work has been devoted to the analysis of second-order non-autonomous evolution equations with nonlocal initial data.
We have established existence and regularity results for both damped and undamped problems by combining fundamental solution techniques, regularity theory for non-autonomous forms, and fixed-point arguments.
The applicability of the developed framework has been illustrated through model problems inspired by realistic physical scenarios, including wave propagation, population dynamics with memory effects, and nonlocal semilinear phenomena.
Overall, the results demonstrate how abstract form methods provide solutions with optimal regularity and offer a systematic approach for treating concrete classes of partial differential equations.


\begin{thebibliography}{99}  
\bibitem{Achache2018}
M. Achache, \textit{Maximal regularity for non-autonomous evolution equations}, PhD Thesis, Université de Bordeaux, 2018.

\bibitem{MA2020}
M. Achache, \textit{Maximal regularity for the damped wave equations. Journal of Elliptic and Parabolic Equations}, \textbf{6}, (2020), 835--870.

\bibitem{AgarwalLeiva2024}
R.P. Agarwal, H. Leiva, L. Riera, S. Lalvay, \textit{Existence of solutions for impulsive neutral semilinear evolution equations with nonlocal conditions}, Discontinuity, Nonlinearity, and Complexity, \textbf{13}, (2024), 333-349. \url{https://doi.org/10.5890/DNC.2024.06.011}

\bibitem{arendt}
W. Arendt, H. Vogt, J. Voig, \textit{Form Methods for Evolution Equations}, Lecture Notes, 2011.

\bibitem{Arendth2014}
W. Arendt, D. Dier, H. Laasri, E.M. Ouhabaz, \textit{Maximal regularity for evolution equations governed by non-autonomous forms}, Advances in Differential Equations, \textbf{19}, (2014) 1043--1066.

\bibitem{Arendth2016}
W. Arendt, S. Monniaux, \textit{Maximal regularity for non-autonomous Robin boundary conditions}, Math. Nachr. \textbf{289} (2016), 1325--1340.

\bibitem{Auscher2016}
P. Auscher, \textit{On non-autonomous maximal regularity for elliptic operators in divergence form}, Arch. Math. \textbf{107} (2016), 271--284.


\bibitem{AHLMT2002} 
P. Auscher, S. Hofmann, M. Lacey, A. McIntosh, P. Tchamitchian, \textit{The solution of the Kato square root problem for second order elliptic operators on $\mathbb{R}^n$}, Ann. Math. \textbf{156} (2002), 633--654.

\bibitem{BezerraSastreSilva2024}
F.D. Bezerra, S. Sastre-Gomez, S.H. da Silva, \textit{Upper semicontinuity for a class of nonlocal evolution equations with Neumann condition}, arXiv preprint arXiv:2409.10065, 2024.

\bibitem{Benedetti22}
I. Benedetti, S. Ciani, \textit{Evolution equations with nonlocal initial conditions and superlinear growth}, J. Differ. Equ. \textbf{318} (2022), 270--297.

\bibitem{Benedetti19}
I. Benedetti, E.M. Rocha, \textit{Existence results for evolution equations with superlinear growth}, Topol. Methods Nonlinear Anal. \textbf{54} (2019), 917--936.

\bibitem{Byszewski1991}
L. Byszewski, \textit{Theorems about the existence and uniqueness of solutions of a semilinear evolution nonlocal Cauchy problem}, J. Math. Anal. Appl. \textbf{162} (1991), 494--505.

\bibitem{BoucherifPrecup2007}
A. Boucherif, R. Precup, \textit{Existence of mild solutions for semilinear evolution equations with multipoint nonlocal conditions}, Nonlinear Anal. \textbf{67} (2007), 3504--3514.


\bibitem{Boyer}
F. Boyer, P. Fabrie, \textit{Mathematical tools for the study of the incompressible Navier-Stokes equations and related models}, Springer, New York, 2013.

\bibitem{Budde}
C. Budde, C. Seifert, \textit{Perturbations of non-autonomous second-order abstract Cauchy problems}, Anal. Math. \textbf{50} (2024), 733--755.


\bibitem{Chill2005}
R. Chill, S. Srivastava, \textit{$L^p$-maximal regularity for second order Cauchy problems}, Math. Z. \textbf{251} (2005), 751--781.

\bibitem{S.Chill}
C.J.K. Batty, R. Chill, S. Srivastava, \textit{Maximal regularity for second order non-autonomous Cauchy problems}, Studia Math. \textbf{189} (2008), 205--223.



\bibitem{Colao2022}
V. Colao, L. Muglia, \textit{Solutions to nonlocal evolution equations governed by non-autonomous forms and demicontinuous nonlinearities}, J. Evol. Equ. \textbf{22} (2022), Article 77.

\bibitem{Dier2014}
D. Dier, \textit{Non-autonomous Cauchy problems governed by forms: maximal regularity and invariance}, Ph.D. Thesis, Universität Ulm, Ulm, 2014.


\bibitem{HPP14}
T. Cardinali, E. De Angelis, \textit{Non-autonomous second order differential inclusions with a stabilizing effect}, Results Math. \textbf{76} (2021), Article 8.


\bibitem{Dier2015}
D. Dier, \textit{Non-autonomous maximal regularity for forms of bounded variation}, J. Math. Anal. Appl. \textbf{425} (2015), 33--54.

\bibitem{Dier2017}
D. Dier, R. Zacher, \textit{Non-autonomous maximal regularity in Hilbert spaces}, J. Evol. Equ. \textbf{17} (2017), 883--907.


\bibitem{Fackler2017}
S.J.L. Fackler, \textit{Lions' problem concerning maximal regularity of equations governed by non-autonomous forms}, Ann. Inst. H. Poincaré Anal. Non Linéaire \textbf{34} (2017), 699--709.

\bibitem{Haak2015}
B.H. Haak, E.M. Ouhabaz, \textit{Maximal regularity for non-autonomous evolution equations}, Math. Ann. \textbf{363} (2015), 1117--1145.


\bibitem{Kozak}
M. Kozak, \textit{A fundamental solution of a second-order differential equation in a Banach space}, Univ. Iagel. Acta Math. \textbf{32} (1995), 275--289.



\bibitem{Leray}
J. Leray, J. Schauder, \textit{Topologie et équations fonctionnelles}, Ann. Sci. Éc. Norm. Supér. \textbf{51} (1934), 45--78.


\bibitem{LinTyniZimmermann2024}
Y.-H. Lin, T. Tyni, P. Zimmermann, \textit{Well-posedness and inverse problems for semilinear nonlocal wave equations}, arXiv preprint arXiv:2402.05877, 2024.


\bibitem{Lions1961}
J.L. Lions, \textit{Équations différentielles opérationnelles et problèmes aux limites}, Springer-Verlag, Berlin, 1961.


\bibitem{Lucchetti}
R. Lucchetti, F. Patrone, \textit{On Nemytskii's operator and its application to the lower semicontinuity of integral functionals}, Indiana Univ. Math. J. \textbf{29} (1980), 703--713.


\bibitem{Ntouyas2005} 
S.K. Ntouyas, \textit{Nonlocal initial and boundary value problems: a survey}, in: Handbook of Differential Equations: Ordinary Differential Equations, Vol. 2, Elsevier, Amsterdam, 2005, 461--557.



\bibitem{NguyenTranVu2024}
D. Nguyen Van, K. Tran Dinh, P. Vu, \textit{Stability analysis for a class of semilinear nonlocal evolution equations}, Electron. J. Qual. Theory Differ. Equ. \textbf{2024} (2024), Article 23.

\bibitem{OD2014}
D. Dier, E.M. Ouhabaz, \textit{Maximal regularity for non-autonomous second order Cauchy problems}, Integral Equ. Oper. Theory \textbf{78} (2014), 427--450.

\bibitem{Ouhabaz2015}
E.M. Ouhabaz, \textit{Maximal regularity for non-autonomous evolution equations governed by forms having less regularity}, Arch. Math. \textbf{105} (2015), 79--91.


\bibitem{PaicuVrabie2010}
A. Paicu, I.I. Vrabie, \textit{A class of nonlinear evolution equations subjected to nonlocal initial conditions}, Nonlinear Anal. \textbf{72} (2010), 4091--4100.

\bibitem{Pazzy}
A. Pazy, \textit{Semigroups of linear operators and applications to partial differential equations}, Springer Science \& Business Media, 2012.

\bibitem{SaniLaasri19}
A. Sani, H. Laasri, \textit{Evolution equations governed by Lipschitz continuous non-autonomous forms}, Czech. Math. J. \textbf{65} (2019), 475--491.

\bibitem{SchmitzWalker2024}
L. S. Schmitz, C. Walker, \textit{Recovering initial states in semilinear parabolic problems from time-averages}, arXiv preprint arXiv:2407.03829, 2024.

\bibitem{Surway}
W. Arendt, D. Dier, S.J.L. Fackler, \textit{Lions’ problem on maximal regularity}, Arch. Math. \textbf{109} (2017), 59--72.

\bibitem{WC2007}
W. Arendt, R. Chill, S. Fornaro, C. Poupaud, \textit{$L^p$-maximal regularity for non-autonomous evolution equations}, J. Differ. Equ. \textbf{237} (2007), 1--26.


\bibitem{XuCaraballoValero2024}
J. Xu, T. Caraballo, J. Valero, \textit{Asymptotic behavior of a semilinear problem in heat conduction with long time memory and non-local diffusion}, arXiv preprint arXiv:2407.17923, 2024.

\bibitem{XuColaoMuglia2021}
H.K. Xu, V. Colao, L. Muglia, \textit{Mild solutions of nonlocal semilinear evolution equations on unbounded intervals via approximation solvability method in reflexive Banach spaces}, J. Math. Anal. Appl. \textbf{498} (2021), Article 124938.


\bibitem{19}
D. Daners, \textit{Heat kernel estimates for operators with boundary conditions}, Math. Nachr. \textbf{217} (2000), 13--41.


\bibitem{22}
O. El-Mennaoui, H. Laasri, \textit{A note on the norm-continuity for evolution families arising from non-autonomous forms}, Semigroup Forum \textbf{100} (2020), 205--230.


\bibitem{35}
A. McIntosh, \textit{On representing closed accretive sesquilinear forms as $(A^{1/2}u, A^{1/2}v)$}, Coll. France Semin. Vol. III, 252--267, 1982.

\bibitem{48}

E. Zeidler, \textit{Nonlinear functional analysis and its applications II/B: linear monotone operators}, Springer, New York, 2013.

\end{thebibliography}
\end{document}